\documentclass[10 pt,leqno]{amsart}
\usepackage{amsfonts}
\usepackage{amsthm}\usepackage{amssymb}
\usepackage{latexsym}
\usepackage{amsmath}
\usepackage{mathrsfs}
\usepackage{pifont}

\theoremstyle{plain}
\newtheorem{thm}{Theorem}[section]

\newtheorem {lem} [thm]{Lemma}
\newtheorem {prop}[thm] {Proposition}
\newtheorem {prob}[thm]{Problem}
\newtheorem{question}[thm]{Question}
\newtheorem{cor}[thm]{Corollary}
\newtheorem{con}[thm]{Conjecture}

\theoremstyle{definition}
\newtheorem {defin}[thm] {Definition}
\newtheorem {rem} [thm]{Remark}
\newtheorem {rems} [thm]{Remarks}
\newtheorem {nots} [thm]{Notation}
\newtheorem {examp} [thm]{Example}

\newcommand{\C} {\Bbb C}
\newcommand{\N}{\Bbb N}

\newcommand{\T}{\Bbb T}
\newcommand{\Z}{\Bbb Z}
\newcommand{\ld}{\ldots}

\newcommand{\Aa} {\mathcal{A}}

\newcommand{\Bb} {\mathcal B}

\newcommand{\Tt} {\mathcal T}
\newcommand{\Oo} {\mathcal O}
\newcommand{\Uu} {\mathcal U}
\newcommand{\Vv} {\mathcal V}

\newcommand{\Pp} {\mathcal P} 
\newcommand{\Kk} {\mathcal K}

\newcommand{\be}{{\bf{1}}}
\newcommand {\Tr} {{\textrm{Tr}}}

\newcommand{\dimnuc}{\textup{dim}_{\mathrm{nuc}}\,}
\newcommand{\dr}{\textup{dr}\,}
\newcommand{\Gan}{\Gamma (n)}
\newcommand{\Gakn}{\Gamma_k (n)}
\newcommand{\Gank}{\Gamma (n^k)}
\newcommand{\subs}{\subseteq}

\newcommand{\la}{\lambda}

\newcommand{\eps}{\varepsilon}
\newcommand{\vphi}{\varphi}
\newcommand{\ot}{\otimes}

\newcommand{\asdim}{\textup{asdim }} 

\subjclass[2000]{46L05, 46L85}

\date{\today}

\thanks{{\it Supported by:} EPSRC First Grant EP/G014019/1}

\begin{document}
\author{Wilhelm Winter}
\author{Joachim Zacharias}

\address{School of Mathematical Sciences,  University of Nottingham,Nottingham, NG7 2RD}
\email{wilhelm.winter@nottingham.ac.uk}

%\address{School of Mathematical Sciences,  University of Nottingham,Nottingham, NG7 2RD}
\email{joachim.zacharias@nottingham.ac.uk}

\keywords{noncommutative covering dimension,  classification of $C^*$-algebras, \indent approximation properties, asymptotic dimension}

\title{The nuclear dimension of $C^{*}$-algebras}

\begin{abstract}
\noindent
We introduce the nuclear dimension of a $C^{*}$-algebra; this is a noncommutative version of topological covering dimension based on a modification of the earlier concept of decomposition rank. Our notion behaves well with respect to inductive limits, tensor products, hereditary subalgebras (hence ideals), quotients, and even extensions. It can be computed for many examples; in particular, it is finite for all UCT Kirchberg algebras. In fact, all classes of nuclear $C^{*}$-algebras which have so far been successfully classified consist of examples with finite nuclear dimension, and it turns out that finite nuclear dimension implies many properties relevant for the classification program.  Surprisingly, the concept is also linked to coarse geometry, since for a discrete metric space of bounded geometry the nuclear dimension of the associated uniform Roe algebra is dominated by the asymptotic dimension of the underlying space.
%
%
%
%We define a generalisation of decomposition rank which on commutative $C^*$-algebras 
%coincides with the usual decompostion rank as defined in \cite{KW}. 
%On purely infinite simple nuclear $C^*$-algebras it  attains values bounded by 3. For 
%discrete metric spaces the weak decomposition rank of the associated uniform Roe 
%algebra is dominated by the asymptotic dimension of the space.
\end{abstract}

\maketitle

\section*{Introduction}   \label{introduction}

\noindent
Recent developments in noncommutative topology suggest that dimension type conditions play a crucial role for the understanding of noncommutative spaces and their applications, cf.\ \cite{EGL:simple_AH}, \cite{Win:localizingEC}, \cite{Yu}, \cite{Connes:compact-metric} and \cite{Connes:NCG}. While in the commutative case the various definitions of covering dimension tend to coincide (at least for sufficiently well-behaved spaces), their generalizations to the noncommutative situation yield vastly different notions, such as stable rank, real rank, or decomposition rank (cf.\ \cite{Rie:sr}, \cite{BP:rr0}, \cite{KW}), each of which has turned out to be highly useful and interesting in its own right. The known applications, e.g.\ to the classification of nuclear $C^{*}$-algebras, are all limited to somewhat special situations -- although, from a philosophical point of view, it should be possible to handle many of these in a unified manner.  There are also notions which have not yet been generalized to the noncommutative setting, such as Gromov's asymptotic dimension (and the latter should clearly be accessible from a noncommutative point of view, as it has already been shown to be closely related to the coarse Baum--Connes conjecture). 

The present paper seeks to remedy this situation. We will propose a notion of noncommutative covering dimension which on the one hand is flexible enough to cover large classes of (nuclear) $C^{*}$-algebras, and which on the other hand is intimately related to many other regularity properties of noncommutative spaces. The concept is linked  to the classification program for nuclear $C^{*}$-algebras, as well as to the theory of dynamical systems and to coarse geometry. We hope that it will contribute to a deeper understanding of the interplay between these fields, but also shed new light on the role of dimension type conditions in other areas of noncommutative geometry.

Our \emph{nuclear dimension} is seemingly only a small variation of the decomposition rank, a notion introduced by Kirchberg and the first named author in \cite{KW} (this in turn was based on earlier concepts introduced in \cite{Win:cpr1} and \cite{Win:cpr2}). The decomposition rank models the dimension type condition in terms of a  decomposition property of noncommutative partitions of unity. Nuclear dimension is defined in a similar manner, only now we add a little more flexibility to the partitions of unity under consideration. The outcome is a notion of integer valued covering dimension for nuclear $C^{*}$-algebras, which still coincides with covering dimension of the spectrum in the commutative case, and which still has nice permanence properties.  But now, the added flexibility in the choice of the partitions of unity makes the theory accessible to much larger classes of $C^{*}$-algebras.

The decompostion rank has turned out to be extremely useful for the classification of stably finite, separable, simple, nuclear $C^{*}$-algebras. In fact, all classes of such $C^{*}$-algebras which by now have been classified consist of ones with finite decomposition rank -- and it seems well possible that separable simple $C^{*}$-algebras with finite decomposition rank are entirely classifiable by their $K$-theory data. An important step in this direction was achieved in \cite{Win:dr-Z-stable}, where it was shown that, for separable, simple, unital $C^{*}$-algebras,   finite decomposition rank implies $\mathcal{Z}$-stability, i.e., all such $C^{*}$-algebras absorb the Jiang--Su algebra $\mathcal{Z}$ tensorially. (The Jiang--Su algebra was introduced in \cite{JiaSu:Z}; see \cite{RorWin:Z-revisited} for alternative characterizations.) The decomposition rank can take finite values only for quasidiagonal $C^{*}$-algebras, so its use beyond  the stably finite case of the classification program will be limited. On the other hand, Kirchberg and Phillips have very successfully classified purely infinite simple $C^{*}$-algebras. Although in their initial approach, topological dimension type conditions do not show up explicitly, these nontheless have turned out to be important both in the simple and the nonsimple case, cf.\ \cite{Kir}, \cite{BlanKir}. We will show that the $C^{*}$-algebras covered by Kirchberg--Phillips classification all have finite nuclear dimension, so that our theory covers large parts of the classification program, both in the stably finite and in the purely infinite case. In fact, one of our motivations is to make progress on a unified approach to the classification problem for nuclear $C^{*}$-algebras, i.e., an approach that does not require genuinely different methods in the finite and the  infinite case. 

We have already mentioned  that, in the simple and unital case, finite decomposition rank implies $\mathcal{Z}$-stability. Using results of Kirchberg, we will be able to derive an infinite version of this statement, namely, that a separable simple $C^{*}$-algebra with finite nuclear dimension and no nontrivial trace is purely infinite, hence absorbs the Cuntz algebra $\mathcal{O}_{\infty}$. As of this moment, we do not know whether simplicity and finite nuclear dimension will imply $\mathcal{Z}$-stability in general; however, there are promising results pointing in this direction, see \cite{NgWin:cfp} (where the corona factorization property is confirmed for simple, unital $C^{*}$-algebras with finite nuclear dimension) and Remark~\ref{rem-dichotomy}; cf.\ also Conjecture~\ref{regularity-conjecture} below.    

A natural touchstone for any kind of invariant for $C^{*}$-algebras will be its behavior with respect to standard constructions, such as direct sums, limits, tensor products, quotients, ideals, or hereditary subalgebras. Decomposition rank and nuclear dimension behave equally well in this respect. There is, however, one exception: since finite decomposition rank implies finiteness, the Toeplitz extension shows that finite decomposition rank does not pass from quotients and ideals to extensions in general -- a problem circumvented by the additional flexibility of nuclear dimension. 

The situation for crossed products is more subtle. At this point, we only have partial results about the topological dimension of crossed products; for example, it is known that the transformation group $C^{*}$-algebra of a minimal diffeomorphism on a compact smooth manifold has finite decomposition rank -- and the proof is extremely technical, cf.\ \cite{QLinPhi:minhom-survey} and \cite{Win:subhomdr}. In \cite{TomsWin:minhom}, Toms and the first named author will show that  the transformation group $C^{*}$-algebra of a minimal homeomorphism on an infinite, compact, finite dimensional, metrizable space has finite nuclear dimension -- and this time, the proof is much simpler and more conceptual. (Even more, the methods introduced in this context are an important step towards completing the classification  of $C^{*}$-algebras associated to uniquely ergodic, minimal homeomorphisms on infinite, compact, finite dimensional, metrizable  spaces, as achieved in \cite{TomsWin:minhom}; see also \cite{StrWin:minhom}.)

Another natural situation to consider is when a group satisfies certain geometric dimension type conditions. Here, we face a genuine problem since the (full or reduced) group $C^{*}$-algebra will in general not be nuclear. However, one might as well look at the so-called uniform Roe algebra; it then turns out that if a discrete group (with word length metric) has finite asymptotic dimension in the sense of Gromov, then its uniform Roe algebra has finite nuclear dimension. This statement can be generalized to discrete metric spaces of bounded geometry. At this point it is an open question how much information about the underlying space the Roe algebra actually contains. It will be interesting to approach this question in our context, i.e., analyze what finite nuclear dimension of the Roe algebra means for the underlying space.   The problem is particularly relevant since Yu (in \cite{Yu}) has shown that a group with finite asymptotic dimension satisfies the coarse Baum--Connes conjecture. By now, we know that the  latter also  holds in more general situations, so one might ask whether finite nuclear dimension of the Roe algebra is a strong enough regularity property to ensure  the coarse Baum--Connes conjecture of the underlying group.

Our paper is organized as follows. In Section~\ref{preliminaries} we recall some facts about order zero maps and completely positive approximations of nuclear $C^{*}$-algebras. In Section~\ref{nuclear-dimension} we introduce our nuclear dimension, compare it to the decomposition rank and derive its  permanence properties with respect to inductive limits, quotients, ideals, extensions and hereditary subalgebras. Section~\ref{almost-order-zero} provides a technical result on the special structure of completely positive approximations realizing nuclear dimension; namely, we prove that the outgoing maps can always be chosen to be approximately order zero.  We compare nuclear dimension to Kirchberg's covering number in Section~\ref{cov}. These observations together with a result of Kirchberg are used in Section~\ref{dichotomy} to  obtain a dichotomy result on sufficiently noncommutative $C^{*}$-algebras with finite nuclear dimension: they either have a nontrivial trace or are purely infinite. In Section~\ref{examples} we collect a number of examples both with finite and with infinite nuclear dimension. This list is extended in Sections~\ref{kirchberg-algebras} and \ref{roe-algebras}, where we show that Kirchberg algebras satisfying the Universal Coefficient Theorem have finite nuclear dimension, and that, for a discrete countable metric space of bounded geometry, the nuclear dimension of the associated uniform Roe algebra is dominated by the asymptotic dimension of the space. We close with a number of open problems and possible future developments in Section~\ref{outlook}.

\newpage

\section{Order zero maps}
\label{preliminaries}

\noindent
In this section we recall some facts about order zero maps. These are c.p.\ maps preserving orthogonality; they are particularly well-behaved, and will serve as building blocks of our noncommutative partitions of unity, similar as in \cite{Win:cpr1}, \cite{Win:cpr2} and \cite{KW}.

\begin{defin}
Let $A$ and $B$ be $C^{*}$-algebras, and $\varphi: A \to B$ a c.p.\ map. We say $\varphi$ has order zero, if, for $a,b \in A_{+}$,
\[
a \perp b \Rightarrow \varphi(a) \perp \varphi(b).
\]  
\end{defin}

The following structure theorem for order zero maps was derived in \cite{WinZac:order-zero} (based on results from \cite{Wol:disjointness}, and generalizing \cite[1.2]{Win:fintopdim}, which only covers the case of finite-dimensional domains).

\begin{thm}
\label{order-zero-structure}
Let $A$ and $B$ be $C^{*}$-algebras and $\varphi:A \to B$ a c.p.\ order zero map. Let $C:= C^{*}(\varphi(A)) \subset B$, then there is a positive element $h \in \mathcal{M}(C) \cap C'$ with $\|h\|= \|\varphi\|$ and a $*$-homomorphism
\[
\pi_{\varphi}:A \to \mathcal{M}(C) \cap \{h\}' \subset B^{**}
\]
such that
\[
\pi_{\varphi}(a)h = \varphi(a) \mbox{ for } a \in A.
\]
If $A$ is unital, then $h = \varphi(\be_{A}) \in C$.
\end{thm}

In the situation of the preceding theorem, we call $\pi_{\varphi}$ the canonical supporting $*$-homo\-morphism of $\varphi$.

We shall have use for the following easy consequence of Theorem~\ref{order-zero-structure}, cf.\ \cite{WinZac:order-zero}.

\begin{cor}
\label{order-zero-traces}
Let $A$, $B$ be $C^{*}$-algebras and $\psi:A \to B$ a c.p.c.\ order zero map. If $\tau$ is a  positive tracial functional on $B$, then $\tau \circ \psi$ is a positive tracial functional on A. 
\end{cor}

By \cite[1.2.3]{Win:cpr2}, order zero maps with finite-dimensional domains can be described in terms of generators and relations which are weakly
stable in the sense of \cite{Lor:lifting}.  The following is a straightforward reformulation of 
\cite[Proposition~2.5]{KW} in this context.

\begin{prop}
\label{order-zero-weakly-stable}
Let $F$ be a finite dimensional $C^{*}$-algebra. For any $\eta>0$ there is $\delta>0$ such that 
the following holds: If $A$ is a $C^{*}$-algebra and $\varphi:F \to A$ a c.p.c.\ order zero map, 
and if $d \in A^{+}$ is a positive contraction in the unitization of $A$ satisfying 
$\|[d, \varphi(x)]\| \le \delta \|x\|$ for all $x \in F$, then there is a c.p.c.\ order 
zero map $\hat{\varphi}: F \to A$ such that 
$\| \hat{\varphi}(x) - d^{\frac{1}{2}} \varphi(x) d^{\frac{1}{2}}\| \le \eta \|x\|$ for all $x \in F$.  
\end{prop}

\section{Nuclear dimension}       
\label{nuclear-dimension}

%Throughout the paper $A$ will denote a $C^*$-algebra, usually nuclear. It is well-known 
%that nuclearity is equivalent to the existence of a net of cp approximations 
%$(\psi_{\la},\vphi_{\la},F_{\la})$, where $\psi_{\la}: A \to F_{\la}$, 
%$ \vphi_{\la}: F_{\la} \to A$ are completely positive contractions such that 
%$\psi_{\la} \circ \vphi_{\la} (a) \to a$ for all $a \in A$, uniformly on finite 
%subsets of $A$ (also called an approximating net). 
%In \cite{KW} Kirchberg and the first named author introduced the decomposition rank 
%defined in terms of such approximating nets. Philosophically, this corresponds to 
%regarding $(\psi_{\la}, F_{\la},\vphi_{\la})$ as a noncommutative partition of unity. 
%In this paper we will consider a weaker version of this concept. In particular we 
%will be interested in cp approximations  $(\psi,\vphi,F)$ and approximating nets  
%$(\psi_{\la},\vphi_{\la},F_{\la})$ where the $\psi$'s and $\vphi$'s are still 
%completely positive but the $\vphi$'s not necessarily contractive.

\noindent
Below we define our notion of noncommutative dimension, compare it to other concepts such as topological covering dimension or decomposition rank, and derive its most important permanence properties.

\begin{defin} \label{weakdecomp}
A  $C^*$-algebra $A$ has nuclear dimension at most $n$, if there exists a net 
$(F_{\la},\psi_{\la},\vphi_{\la})_{\lambda \in \Lambda}$ such that the $F_{\lambda}$ are finite-dimensional $C^{*}$-algebras, and such that $\psi_{\la}: A \to F_{\la}$ and  
$\vphi_{\la}: F_{\la} \to A$ are completely positive maps satisfying
\begin{enumerate}
\item $\psi_{\la} \circ \vphi_{\la} (a) \to a$ uniformly on finite subsets of $A$;
\item $\| \psi_{\la} \| \leq 1$; %$\| \phi_{\la} \| \leq n+1$ ; 
\item for each $\lambda$,  $F_{\lambda}$ decomposes into $n+1$ ideals $F_{\la}= F_{\la}^{(0)} \oplus \ldots \oplus F_{\la}^{(n)}$ 
such that $\vphi_{\la}|_{F_{\la}^{(i)}} $ is a c.p.c.\ order zero map for $i=0,1, \ldots ,n$.
\end{enumerate}
We write $\dimnuc A\leq n$ in this case and refer to the maps 
$\varphi_{\lambda}$ as  piecewise contractive $n$-decomposable c.p.\ maps, and to the triples 
$(F_{\la}, \psi_{\la}, \vphi_{\la})$ as piecewise contractive $n$-decomposable 
c.p.\ approximations.
\end{defin}

%Note that for $(F_{\la},\psi_{\la},\vphi_{\la})$ as in this definition we must have 
%$\|\vphi_{\la}\| \leq n+1$. If we add the requirement $\|\vphi_{\la}\|\leq 1$ then 
%we obtain the ordinary decomposition rank $ \textup{dr}\,(A)$ (\cite{KW}). This seemingly 
%small variation has significant consequences.
%The following properties can be verified immediately. 

\begin{rems} \label{basic}
Let $A$ be a $C^{*}$-algebra. 
\begin{enumerate}
\item If $\dimnuc A \le n < \infty$, there is a system $(F_{\lambda},\psi_{\lambda},\varphi_{\lambda})_{\Lambda}$ of piecewise contractive $n$-decomposable c.p.\ approximations for $A$; since each map  $\varphi_{\lambda}$ is a sum of at most $n+1$ c.p.c.\ maps, and since the $\psi_{\lambda}$ are c.p.c.\ maps, the norms of the compositions $\varphi_{\lambda} \psi_{\lambda}$ are uniformly bounded by $n+1$. It is straightforward to check that this implies that $A$ has the completely positive approximation property, so that $A$ is nuclear. 
\item Recall from \cite{KW} that the decomposition rank is defined almost exactly as the nuclear dimension, with only the -- seemingly small -- extra condition that the maps $\varphi_{\lambda}$ themselves are contractive. It is therefore trivial that $\dimnuc A \leq \dr A$.
\item It is also trivial that $\dimnuc A =0$ iff $\dr A=0$. Moreover, this happens iff $A$ is an AF algebra (cf.\ \cite[Example~4.1]{KW}).
\item We will see later that nuclear dimension and decomposition rank in general do not coincide, so that we cannot generally choose the maps $\psi_{\lambda}$ and $\varphi_{\lambda}$ in the approximations $(F_{\lambda},\psi_{\lambda},\varphi_{\lambda})_{\Lambda}$ of \ref{weakdecomp} to be contractive. We can, however, always modify the approximations such that the compositions $\varphi_{\lambda} \psi_{\lambda}$ are indeed contractions. 
If $A$ is unital, then $\varphi_{\lambda} \psi_{\lambda}(\be_{A}) \to \be_{A}$, so it will suffice to replace $\varphi_{\lambda}$ by $\|\varphi_{\lambda}\psi_{\lambda}(\be_{A})\|^{-1} \cdot \varphi_{\lambda}$.    In the nonunital case choose an approximate unit $(u_{\sigma})_{\Sigma}$ for $A$, and replace the net $\Lambda$ by the double-indexed net $\Lambda \times \Sigma$, the maps $\psi_{\lambda}$ by $\psi_{\lambda,\sigma}:= \psi_{\lambda}(u_{\sigma}^{\frac{1}{2}}\, . \,u_{\sigma}^{\frac{1}{2}})$ and the maps $\varphi_{\lambda}$ by $\varphi_{\lambda,\sigma}:= \|\varphi_{\lambda}\psi_{\lambda}(u_{\sigma})\|^{-1} \cdot \varphi_{\lambda}$. It is then straightforward to check that these new approximations still form a system of piecewise contractive, $n$-decomposable c.p.\ approximations with the additional property that the compositions $\varphi_{\lambda,\sigma} \psi_{\lambda,\sigma}$ are contractions.
\item Note that we did not ask $A$ to be separable in  Definition~\ref{weakdecomp}. While \cite[Definition~3.1]{KW} was formulated only for separable $C^{*}$-algebras, it clearly makes sense in the general situation as well; moreover, several of the basic results of \cite{KW} still hold in the nonseparable case. In the present paper, we did not want to make any restrictions along these lines, since some of our main examples will be nonseparable (cf.\ Section~\ref{roe-algebras}). However,  for many applications one can nontheless restrict to the separable case, cf.\ Proposition~\ref{separable-wdr} below. 
\end{enumerate}
\end{rems}

The following permanence properties are derived just as for the completely positive rank or for the decomposition rank, cf.\ \cite[Section~3]{Win:cpr1} and  \cite[Section~3]{KW}. Note that there is no need to specify the tensor product in \ref{permanence}(ii), since the values can be finite only for nuclear $C^{*}$-algebras.

\begin{prop}
\label{permanence}
Let $A$, $B$, $C$, $D$ and $E$  be $C^{*}$-algebras; suppose $C=\lim_{\to} C_{i}$ is an inductive limit of $C^{*}$-algebras and $D$ is a quotient of $E$. Then,
\begin{enumerate}
\item  $\dimnuc(A \oplus B) = \max(\dimnuc A, \dimnuc B)$ 
\item  $\dimnuc (A \otimes B ) \leq (\dimnuc A +1)(\dimnuc B +1) -1$; if $B$ is an AF algebra, then $\dimnuc (A \otimes B) \le \dimnuc A$
\item  $\dimnuc C \leq \liminf (\dimnuc C_i)$
\item $\dimnuc D \le \dimnuc E$.
\end{enumerate}
\end{prop}

%\begin{rem}
%\label{invertible-approximation}
%In Definition \ref{weakdecomp}, if $A$ is unital we may replace the approximations $(F_{\lambda},\psi_{\lambda},\varphi_{\lambda})$ by  $(p_{\lambda}F_{\lambda}p_{\lambda},p_{\lambda}\psi_{\lambda}(\, . \,)p_{\lambda},\varphi_{\lambda}|_{p_{\lambda}F_{\lambda}p_{\lambda}})$, the latter still being an piecewise contractive $n$-decomposable c.p.\ approximation. This way, we may always assume the elements $\psi_{\lambda}(\be_{A})$ to be invertible in $F_{\lambda}$. 
%\end{rem}

Just as the decomposition rank, nuclear dimension agrees with covering dimension of the spectrum in the separable commutative case. In the nonseparable case, nuclear dimension and decomposition rank still coincide, and they agree with the respective definition of covering dimension. The only reason why we distinguish between the separable and the nonseparable case is that the various characterizations of dimension tend to disagree for  spaces which are not second countable.

\begin{prop} \label{commutative}
Let $X$ be a locally compact Hausdorff space. Then,  
\[
\dimnuc C_{0}(X) =\dr C_{0}(X).
\]
In particular, if $X$ is second countable, we have  
\[
\dimnuc C_{0}(X) =\dr C_{0}(X) = \dim X.
\]
\end{prop}

\begin{proof}
We have $\dimnuc C_{0}(X) \le \dr C_{0}(X)$ by Remark~\ref{basic}. For the reverse estimate, 
let us assume that $\dimnuc C_{0}(X) = n < \infty$. Suppose $\mathcal{F} \subset C_{0}(X)$ 
is a finite subset of positive normalized elements, and that $\varepsilon>0$ is 
given. We may assume that the elements of $\mathcal{F}$ have compact support and that  there is a positive normalized function $h \in C_{0}(X)$ such that $h a = a$ for all $a \in \mathcal{F}$.

Choose a piecewise contractive $n$-decomposable c.p.\ approximation 
$(F= F^{(0)}\oplus \ldots \oplus F^{(n)},\psi,\varphi)$ for $\mathcal{F} \cup \{h\}$ 
within $\varepsilon/2$. Since $\varphi$ has order zero on each matrix block 
of $F$, we see from \cite[Remark~2.16(ii)]{Win:cpr1} that $F$ is commutative. By cutting down 
$F$ to the hereditary subalgebra generated by $\psi(h)$, we may assume 
that $\psi(h)$ is invertible in $F$. Define  c.p.\ maps 
\[
\hat{\psi}: C_{0}(X) \to F \mbox{ and } \hat{\varphi}: F \to C_{0}(X)
\] 
by 
\[
\hat{\psi}(f) := \psi(h)^{-\frac{1}{2}} \psi(h f) \psi(h)^{-\frac{1}{2}} \mbox{ for } f \in C_{0}(X)
\]
and
\[
\hat{\varphi}(x):= \left(1 - \frac{\varepsilon}{2} \right) \cdot \varphi(\psi(h)^{\frac{1}{2}}  x  \psi(h)^{\frac{1}{2}}) \mbox{ for } x \in F.
\] 
It is clear that $\hat{\psi}$ is  contractive, and that  
$\hat{\varphi}$ is $n$-decomposable with respect to $F= F^{(0)} \oplus \ldots \oplus F^{(n)}$. Moreover, 
\begin{eqnarray*}
\hat{\varphi}(\be_{F}) & = & \hat{\varphi}\hat{\psi}(h) \\
& = & \left(1 - \frac{\varepsilon}{2} \right) \cdot \varphi \psi (h) \\
& \le & \left(1 - \frac{\varepsilon}{2} \right) \left(1 +  \frac{\varepsilon}{2} \right) \cdot  h \\
& \le & \be ,
\end{eqnarray*}
whence $\hat{\varphi}$ is contractive. Finally, we have
\begin{eqnarray*}
\|\hat{\varphi} \hat{\psi} (f) - f\| & \le &  \|\hat{\varphi} \hat{\psi} (f) - \varphi \psi(f) \| + \|\varphi \psi(f)- f\| \\
& < & \left\|\left(1 - \frac{\varepsilon}{2}\right) \cdot \varphi \psi(hf)  - \varphi \psi(f) \right\| + \frac{\varepsilon}{2} \\
& \le & \varepsilon
\end{eqnarray*}
for $f \in \mathcal{F}$, so $(F,\hat{\psi},\hat{\varphi})$ is an $n$-decomposable 
c.p.c.\ approximation for $\mathcal{F}$ within $\varepsilon$. Therefore,  $\dr C_{0}(X) \le n$.

The statement about the second countable case is \cite[Proposition~3.3]{KW}.
\end{proof}

As we pointed out in Remark~\ref{basic}, we do not wish to impose any separability restrictions on our definition of nuclear dimension. However, in many situations one can nontheless restrict to the separable case, using the following observation.

We next show that, just like for decomposition rank, finite nuclear dimension passes to hereditary subalgebras.  Combined with Brown's Theorem, this result shows that nuclear dimension is a stable invariant, cf.\ Corollary~\ref{her-cor} below.

\begin{prop} 
\label{hereditary}
$\dimnuc B \leq \dimnuc A$ when $B \subs A$ is a hereditary $C^{*}$-subalgebra.
\end{prop}

\begin{proof}
We may assume $n:= \dimnuc A$ to be finite, for otherwise there is nothing to show. 
Let $b_{1}, \ldots, b_{m} \in B_{+}$ be normalized elements and let $\varepsilon>0$ 
be given. We have to find a piecewise contractive $n$-decomposable c.p.\ approximation 
(of $B$) for $\{b_{1},\ldots,b_{m}\}$ within $\varepsilon$. 

Using an idempotent approximate unit, by slightly perturbing the $b_{j}$ we may 
(as in \cite[Remark~3.2(ii)]{KW}) assume that there are positive normalized 
elements $h_{0}, \,h_{1} \in B_{+}$ such that 
\[
h_{0}h_{1} = h_{1} \mbox{ and } h_{1}b_{j} = b_{j}
\]
for $j=1, \ldots,m$. 

Set 
\[
\eta:= \min \left\{\frac{\varepsilon^{8}}{13 (n+1)}, \, \frac{1}{2^{16}} \right\}
\]
and choose a piecewise contractive $n$-decomposable c.p.\ approximation 
\[
(F=F^{(0)} \oplus \ldots F^{(n)},\psi,\varphi)
\]
(of $A$) for $\{h_{0},h_{1}, b_{1}, \ldots, b_{m}\}$ within $\eta$. 

Define a projection $p \in F$ by 
\[
p:= g_{\eta^{\frac{1}{2}}}(\psi(h_{1})),
\] 
where $g_{\eta^{\frac{1}{2}}}$ is given by
\[
g _{\eta^{\frac{1}{2}}}(t):= 
\left\{ \begin{array}{ll}
0& \mbox{for } t < \eta^{\frac{1}{2}}, \\
1 & \mbox{for } t \ge \eta^{\frac{1}{2}}. \end{array} \right.
\]
Set 
\[
\hat{F}:= pFp, \,  \hat{F}^{(i)}:= pF^{(i)}p \mbox{ and } p^{(i)}:= \be_{F^{(i)}}p
\]
for $i \in \{0, \ldots, n\}$ and define a c.p.c.\ map 
\[
\hat{\psi}:B \to \hat{F}
\] 
by 
\[
\hat{\psi}(b):= p \psi(b)p,  \, b \in B.
\] 
For $i \in \{0, \ldots, n\}$ we have
\begin{eqnarray*}
\|\varphi^{(i)}(p^{(i)}) (\be - h_{0})\|  & = & \|(\be - h_{0}) \varphi^{(i)}(p^{(i)})^{2} (\be - h_{0})\|^{\frac{1}{2}} \\
& \le & \|(\be - h_{0}) \varphi^{(i)}(p^{(i)}) (\be - h_{0})\|^{\frac{1}{2}} \\
& \le & \|(\be - h_{0}) \varphi(p) (\be - h_{0})\|^{\frac{1}{2}} \\
& \le & \left(\frac{1}{\eta^{\frac{1}{2}}} \| (\be - h_{0}) \varphi \psi(h_{1}) (\be - h_{0}) \| \right)^{\frac{1}{2}} \\
& \le & \left(\frac{\eta}{\eta^{\frac{1}{2}}} \right)^{\frac{1}{2}} \\
& = & \eta^{\frac{1}{4}} \left(\le \frac{1}{16}\right).
\end{eqnarray*}
Now by \cite[Lemma~3.6]{KW} (applied to $\varphi^{(i)}|_{\hat{F}^{(i)}}$ in 
place of $\varphi$ and $h_{0}$ in place of $h$) there are c.p.c.\ order zero maps 
\[
\hat{\varphi}^{(i)}: \hat{F}^{(i)} \to \overline{h_{0}Ah_{0}} \subset B
\]
such that 
\[
\| \hat{\varphi}^{(i)}(x) - \varphi^{(i)}(x)\| \le 8 \eta^{\frac{1}{8}} \|x\|
\]
for all $0 \le x \in \hat{F}^{(i)}$ and $i \in \{0, \ldots, n\}$. Set 
\[
\hat{\varphi}:= \sum_{i=0}^{n} \hat{\varphi}^{(i)}: \hat{F} \to B.
\] 
The map $\hat{\varphi}$ is a sum of $n+1$ c.p.c.\ order zero maps by 
construction, and we have 
\begin{equation}
\|\hat{\varphi} \hat{\psi}(b_{j}) - \varphi \hat{\psi}(b_{j})\| \le 8 (n+1) \eta^{\frac{1}{8}}, \, j = 1, \ldots, m. \label{her1}
\end{equation} 
To check that $\hat{\varphi} \hat{\psi}(b_{j})$ is close to $b_{j}$, note 
first that
\begin{eqnarray*}
\| \varphi((\be_{F}-p) \psi(b_{j})) \| & \le & \|\varphi((\be_{F}-p) \psi(b_{j})) \varphi(\psi(b_{j})(\be_{F}-p))\|^{\frac{1}{2}} \\
& \le & \|\varphi((\be_{F}-p) \psi(b_{j})^{2}(\be_{F}-p))\|^{\frac{1}{2}} \\
& \le & \|\varphi((\be_{F}-p) \psi(h_{1})(\be_{F}-p))\|^{\frac{1}{2}} \\
& \le & ((n+1) \eta^{\frac{1}{2}})^{\frac{1}{2}} \\
& \le & (n+1) \eta^{\frac{1}{4}}
\end{eqnarray*}
for $j = 1, \ldots, m$, which in particular implies that 
\[
\|\varphi([p, \psi(b_{j})]) \| \le 2 (n+1) \eta^{\frac{1}{4}}.
\] 
We now obtain 
\begin{eqnarray}
\|\varphi \psi(b_{j}) - \varphi \hat{\psi}(b_{j}) \| & \le & \| \varphi(\psi(b_{j}) - p \psi(b_{j}) + \psi(b_{j})p - p \psi(b_{j})p ) \| + 2 (n+1) \eta^{\frac{1}{4}} \nonumber \\
& \le & 4 (n+1) \eta^{\frac{1}{4}} \label{her2}
\end{eqnarray}
for $j = 1, \ldots,m$, whence 
\begin{eqnarray*}
\| \hat{\varphi} \hat{\psi} (b_{j}) - b_{j} \| & \le & \| \hat{\varphi} \hat{\psi} (b_{j}) - \varphi \hat{\psi} (b_{j}) \|  +  \| \varphi \hat{\psi} (b_{j}) - \varphi \psi (b_{j}) \| + \| \varphi \psi (b_{j}) - b_{j} \| \\
& \stackrel{\eqref{her1}, \eqref{her2}}{<} & 8 (n+1) \eta^{\frac{1}{8}} + 4 (n+1) \eta^{\frac{1}{4}} + \eta \\
& < & \varepsilon.
\end{eqnarray*}
Therefore, the approximation $(\hat{F},\hat{\psi}, \hat{\varphi})$ is as desired.
\end{proof}

\begin{prop}
\label{separable-wdr}
Let $B$ be a $C^{*}$-algebra. For any countable subset $S \subset B$ there is a 
separable $C^{*}$-subalgebra $C \subset B$ such that $S \subset C$ and $\dimnuc C \le \dimnuc B$. 
\end{prop}

\begin{proof}
Let $\dimnuc B =n < \infty$. Set $S_{0}:= S$ and choose 
\[
(F_{0,\lambda},\psi_{0,\lambda},\varphi_{0,\lambda})_{\lambda \in \mathbb{N}},
\]
a system of piecewise contractive $n$-decomposable c.p.\ approximations (of $B$) for $S_{0}$. 

If $S_{k} \subset B$ and 
\[
(F_{k,\lambda},\psi_{k,\lambda},\varphi_{k,\lambda})_{\lambda \in \mathbb{N}}
\]
have been constructed, choose a countable dense subset
\[
S_{k+1} \subset C^{*}\left(\bigcup_{l \le k, \,  \lambda \in \mathbb{N}} \varphi_{l,\lambda}(F_{l,\lambda}) \cup S_{k}\right) \subset B
\]
and choose 
\[
(F_{k+1,\lambda},\psi_{k+1,\lambda},\varphi_{k+1,\lambda})_{\lambda \in \mathbb{N}},
\]
a system of piecewise contractive, $n$-decomposable c.p.\ approximations (of $B$) for $S_{k+1}$. Continue inductively and define
\[
C:= \overline{\bigcup_{k \in \mathbb{N}}S_{k}};
\]
it is straightforward to check that $C$ has the right properties, and that a system of piecewise contractive, $n$-decomposable c.p.\ approximations of $C$ is given by 
\[
(F_{k,\lambda},\psi_{k,\lambda},\varphi_{k,\lambda})_{k,\lambda \in \mathbb{N}}.
\]
 \end{proof}

\begin{rem}
\label{separable-wdr-remark}
A small modification of the proof above even shows the following:

Let $A$ be a $C^{*}$-algebra and $B \subset A$ a hereditary $C^{*}$-subalgebra. 
For any countable subset $S \subset A$ there is a separable $C^{*}$-subalgebra 
$D \subset A$ such that $S \subset D$ and such that $C:= D \cap B$ (which is 
hereditary in $D$) satisfies $\dimnuc C \le \dimnuc B$. 

If, additionally, $B$ is full in $A$, then $C$ may be taken to be full in $D$.
\end{rem}

\begin{cor}
\label{her-cor}
Let $A$ be a $C^{*}$-algebra.
\begin{itemize}
\item[(i)] For any $r \in \mathbb{N}$ we have $\dimnuc A = \dimnuc (M_{r} \otimes A) = \dimnuc ( \mathcal{K} \otimes A)$. 
\item[(ii)] If $B \subset A$ is a full hereditary $C^{*}$-subalgebra, then $\dimnuc B = \dimnuc A$. 
\end{itemize}
\end{cor}

\begin{proof}
(i) We have $\dimnuc A \le \dimnuc (M_{r} \otimes A) \le \dimnuc(\mathcal{K} \otimes A)$ 
by Proposition~\ref{hereditary}  and $\dimnuc(\mathcal{K} \otimes A) \le \dimnuc A$ by Proposition~\ref{permanence}.

(ii) We have $n:= \dimnuc B \le \dimnuc A$ by Proposition~\ref{hereditary}, so it 
remains to show that $\dimnuc A \le \dimnuc B$. 

Given $a_{1}, \ldots, a_{m} \in A_{+}$, by Remark~\ref{separable-wdr-remark} 
there is a separable $C^{*}$-subalgebra $D \subset A$ such that $\{a_{1}, \ldots, a_{m}\} \subset D$, 
such that $C := D \cap B$ is full in $D$ and such that $\dimnuc C \le \dimnuc B$.  

Now by Brown's Theorem \cite[Theorem~2.8]{Br}, we have $\mathcal{K} \otimes C \cong \mathcal{K} \otimes D$, 
hence $\dimnuc D = \dimnuc C (\le \dimnuc B)$ by part (i) of the corollary. We may thus find 
arbitrarily close piecewise contractive $n$-decomposable c.p.\ approximations for $a_{1}, \ldots, a_{m}$. 
\end{proof}

%The next result describes the behaviour of the decomposition rank under quotients and extensions. 
%Contrary to the case of the ordinary decomposition rank we do not have to impose any restrictions on the 
%extensions involved in part (ii).

We are now ready to describe the first significant difference between decomposition rank and nuclear dimension. We already know that both theories behave well with respect to quotients and ideals;  it has been observed in \cite{KW} that finite decomposition rank passes to quasidiagonal extensions, and that one cannot expect a general statement in this context. The additional flexibility in the definition of nuclear dimension, however, ensures that finite nuclear dimension indeed passes to arbitrary extensions. So we obtain a  noncommutative version of the sum theorem for covering dimension, cf.\ \cite[III.2.B)]{HurWal:dimension}. This behavior will also make large new classes of $C^{*}$-algebras accessible to our theory, cf.\ Example~\ref{toeplitz-example} and Sections~\ref{kirchberg-algebras} and \ref{roe-algebras} below.

\begin{prop} \label{quotext}
Let $0 \to J \to E \to A \to 0$ be an exact sequence of $C^*$-algebras. Then, 
\[
\max \{ \dimnuc A , \, \dimnuc J \} \leq  \dimnuc E \leq \dimnuc A + \dimnuc J +1.
\]
\end{prop}

\begin{proof}
The first inequality follows from  Propositions~\ref{permanence} and \ref{hereditary}.

For the second inequality, we may assume that both $m:= \dimnuc J$ and $n:= \dimnuc A$ are finite, for otherwise 
there is nothing to show. Let positive and normalized elements $e_{1}, \ldots,e_{k} \in E$ 
and $\varepsilon >0$ be given. 

Choose a piecewise contractive $n$-decomposable c.p.\ approximation 
\[
(F_{A}=F_{A}^{(0)}\oplus \ldots \oplus F_{A}^{(n)}, \psi_{A}, \varphi_{A})
\]
(of $A$) for $\{\pi(e_{1}), \ldots,\pi(e_{k})\}$ within $\frac{\varepsilon}{5}$. 
By \cite[Proposition~1.2.4]{Win:cpr2} (essentially using that cones over finite-dimensional $C^{*}$-algebras are projective), each $\varphi_{A}^{(j)}$ lifts to a c.p.c.\ order zero 
map 
\[
\bar{\varphi}_{A}^{(j)}: F_{A}^{(j)} \to E,
\] 
so that 
\[
\bar{\varphi}_{A}:= \sum_{j=0}^{n} \bar{\varphi}_{A}^{(j)}
\] 
will be a piecewise contractive $n$-decomposable c.p.\ lift of $\varphi_{A}$.  

From \cite[1.2.3]{Win:cpr2}, we know that the relations defining order zero maps are weakly stable; 
this in particular implies that there is $\delta>0$ such that the assertion of 
Proposition~\ref{order-zero-weakly-stable} holds for each $F_{A}^{(j)}$ in place 
of $F$ and $\frac{\varepsilon}{5(n+1)}$ in place of $\eta$. 

Using a quasicentral approximate unit for $J$ relative to $E$, it is straightforward 
to find   a positive normalized element $h \in J$ such that the following hold:
\begin{itemize}
\item[(a)] $\|[(\be - h), \bar{\varphi}^{(j)}_{A}(x)]\| \le \delta \|x\|$ for $x \in F_{A}^{(j)}$, $j=0, \ldots,n$
\item[(b)] $\|h^{\frac{1}{2}} e_{l} h^{\frac{1}{2}} + (\be - h)^{\frac{1}{2}} e_{l} (\be - h)^{\frac{1}{2}} - e_{l} \| < \frac{\varepsilon}{5} $ for $l=1, \ldots, k$
\item[(c)] $\|  (\be - h)^{\frac{1}{2}} (\bar{\varphi}_{A}\psi_{A} \pi(e_{l})   -  e_{l}) (\be - h)^{\frac{1}{2}}\| < \frac{2 \varepsilon}{5}$ for $l=1, \ldots, k$.
\end{itemize}
(To obtain (c), we use that 
\[
\|\pi(\bar{\varphi}_{A} \psi_{A} \pi(e_{l})- e_{l}) \| = \|\varphi_{A} \psi_{A} \pi(e_{l}) - e_{l} \| < \frac{\varepsilon}{5},
\]
whence $\bar{\varphi}_{A} \psi_{A} \pi(e_{l}) - e_{l}$ is at most $\frac{\varepsilon}{5}$ away from $J$.) 

Now by (a) and Proposition~\ref{order-zero-weakly-stable} there are c.p.c.\ order zero maps 
\[
\hat{\varphi}_{A}^{(j)}: F_{A}^{(j)} \to E
\]
such that
\[
\|\hat{\varphi}^{(j)}_{A}(x) - (\be-h)^{\frac{1}{2}} \bar{\varphi}_{A}^{(j)}(x) (\be-h)^{\frac{1}{2}} \| \le \frac{\varepsilon}{5(n+1)} \|x\|
\]
for $x \in F_{A}^{(j)}$, $j= 0, \ldots,n$; set
\[
\hat{\varphi}_{A}:= \sum_{j=0}^{n} \hat{\varphi}_{A}^{(j)},
\]
then
\[
\|\hat{\varphi}_{A}(x) - (\be -h)^{\frac{1}{2}} \bar{\varphi}_{A}(x) (\be - h)^{\frac{1}{2}}\| \le \frac{\varepsilon}{5} \|x\| \mbox{ for } x \in F_{A}.
\]
Next, choose a piecewise contractive $m$-decomposable c.p.\ approximation 
\[
(F_{J}=F_{J}^{(0)} \oplus \ldots \oplus F_{J}^{(m)}, \psi_{J},\varphi_{J})
\] 
(of $J$) for $\{h^{\frac{1}{2}}e_{l} h^{\frac{1}{2}} \mid l=1, \ldots,k\}$ 
within $\frac{\varepsilon}{5}$. 

Set
\[
F:= F_{J} \oplus F_{A} , \, \psi(\, . \,):= \psi_{J}(h^{\frac{1}{2}}\, . \,h^{\frac{1}{2}}) \oplus \psi_{A}\pi(\, . \,) \mbox{ and } \varphi:= \varphi_{J} + \hat{\varphi}_{A},
\]
then $\psi$ is c.p.c.\ and $\varphi$ is piecewise contractive c.p.; $\varphi$ is 
$(m+n+1)$-decomposable with respect to  $F= \bigoplus_{j=0}^{m+n+1} F^{(j)}$, where
\[
F^{(j)}:= \left\{
\begin{array}{ll}
F_{J}^{(j)} & \mbox{for } j = 0, \ldots,m \\
F_{A}^{(j-m-1)} & \mbox{for } j=m+1, \ldots, m+n+1.
\end{array} 
\right.
\]
It remains to be checked that $(F,\psi,\varphi)$ indeed approximates the 
$e_{l}$ within $\varepsilon$, i.e., 
\begin{eqnarray*}
\|\varphi \psi(e_{l}) - e_{l}\|& \stackrel{\mathrm{(b)}}{<}& \| \varphi_{J} \psi_{J} (h^{\frac{1}{2}} e_{l} h^{\frac{1}{2}}) - h^{\frac{1}{2}} e_{l} h^{\frac{1}{2}}  \| \\
&& + \| \hat{\varphi}_{A} \psi_{A} \pi(e_{l}) - (\be - h)^{\frac{1}{2}} e_{l} (\be - h)^{\frac{1}{2}} \| \\
&& + \frac{\varepsilon}{5} \\
& \le &  \frac{\varepsilon}{5} + \|  (\be - h)^{\frac{1}{2}}( \bar{\varphi}_{A} \psi_{A} \pi(e_{l}) - e_{l}  )  (\be - h)^{\frac{1}{2}}\| + \frac{\varepsilon}{5} + \frac{\varepsilon}{5} \\
& \stackrel{\mathrm{(c)}}{<} & \varepsilon.
\end{eqnarray*}
\end{proof}

\begin{cor}
Let $A$ be a separable continuous trace $C^{*}$-algebra. Then, 
\[
\dimnuc A = \dr A = \dim \hat{A}.
\]
\end{cor}

\begin{proof}
The proof follows that of \cite[Corollary~3.10]{KW} almost verbatim.
\end{proof}

\begin{rem}
Applying the previous result  to  the minimal unitization $A\tilde{ }$ of a $C^{*}$-algebra $A$, one obtains that 
\[
\dimnuc A\tilde{ } \le \dimnuc A +1.
\]
However, following the lines of \cite[Proposition~3.11]{KW}, one can even show that the nuclear dimension of a $C^{*}$-algebra agrees with that of its smallest unitization.  In the separable commutative case, the respective statement also holds for the maximal compactification. One cannot quite expect a noncommutative generalization of the latter result to  our context, since multiplier algebras in general are not nuclear.
\end{rem}

\section{Almost order zero approximations}
\label{almost-order-zero}

\noindent
In \cite{KW} it was shown that $C^{*}$-algebras with finite decomposition rank are quasidiagonal. The reason was that the $n$-decomposable c.p.c.\ approximations may always be chosen so that the maps $\psi_{\lambda}:A \to F_{\lambda}$ are almost multiplicative, cf.\ \cite[Proposition~5.1]{KW}. In this section, we prove an analogous result for nuclear dimension and piecewise contractive $n$-decomposable c.p.\ approximations, saying that the latter may always be chosen to be almost orthogonality preserving.  We first need a simple technical observation.

\begin{prop}
\label{product-domination}
Let $A$ be a $C^{*}$-algebra, and let $0 \le a \le b$ and $0 \le a' \le b'$ be positive 
elements of norm at most one. Then, $\|aa'\|^{2} \le \|bb'\|$.
\end{prop}

\begin{proof}
We simply estimate
\begin{eqnarray*}
\|bb'\| & \ge & \|b^{\frac{1}{2}}b'bb'b^{\frac{1}{2}}\| \\
& = & \|b^{\frac{1}{2}}b' b^{\frac{1}{2}}\|^{2} \\
& \ge & \|b^{\frac{1}{2}}a' b^{\frac{1}{2}}\|^{2} \\
& = & \| (a')^{\frac{1}{2}} b (a')^{\frac{1}{2}}\|^{2} \\
& \ge &  \| (a')^{\frac{1}{2}} a (a')^{\frac{1}{2}}\|^{2} \\
& \ge & \|a' a a'\|^{2} \\
& \ge & \|a' a^{2} a' \|^{2} \\
& = & \|aa'\|^{2}.
\end{eqnarray*}
\end{proof}

\begin{prop}
\label{almost-order-zero-approximation}
Let $A$ be a $C^{*}$-algebra with $\dimnuc A=n <\infty$. Then, there is a 
system $(F_{\lambda},\psi_{\lambda},\varphi_{\lambda})_{\lambda \in \Lambda}$ of almost 
contractive $n$-decomposable c.p.\ approximations such that the map
\[
\textstyle
\bar{\psi}: A \to \prod_{\Lambda} F_{\lambda} / \bigoplus_{\Lambda} F_{\lambda} 
\]
induced by the $\psi_{\lambda}$ has order zero.
\end{prop}

\begin{proof}
Let us first assume $A$ to be separable. In this case, it will suffice to show the 
following: For any  $0<\varepsilon< \frac{1}{(n+2)^{4}}$ and 
any finite subset $\mathcal{F} \subset A$ of positive normalized 
elements,  there is a piecewise contractive $n$-decomposable c.p.\ approximation 
$(F,\psi,\varphi)$ of $A$ such that 
\[
\| \varphi \psi  (b) - b\| < \varepsilon^{\frac{1}{16}} \mbox{ for } b \in \mathcal{F}
\] 
and
\[
\|\psi (c) \psi(c') \| < \varepsilon^{\frac{1}{16}}
\]
whenever $c,c' \in \mathcal{F}$ satisfy $\|c c'\| < \varepsilon$.

So, let $\varepsilon$ and $\mathcal{F}$ as above be given. Choose a 
piecewise contractive $n$-decomposable c.p.\ approximation 
$(\tilde{F},\tilde{\psi},\tilde{\varphi})$ of $A$ such that 
\[
\| \tilde{\varphi} \tilde{\psi}  (b) - b\| < \varepsilon \mbox{ for } b \in \mathcal{F}.
\]
Write $\tilde{F}= M_{r_{1}} \oplus \ldots \oplus M_{r_{s}}$ and 
denote the respective components of $\tilde{\varphi}$ and $\tilde{\psi}$ 
by $\tilde{\varphi}_{j}$ and $\tilde{\psi}_{j}$, respectively. Define
\begin{eqnarray*}
I & := & \{j \in \{ 1, \ldots, s \} \mid \\
&&\| \tilde{\psi}_{j}  (c) \tilde{\psi}_{j} (c') \| \ge \varepsilon^{- \frac{1}{8}} \| \tilde{\varphi} \tilde{\psi}(c) \tilde{\varphi} \tilde{\psi} (c') \|^{\frac{1}{4}}  \\
& & \mbox{for some } c,c' \in \mathcal{F} \mbox{ with } \|cc'\| < \varepsilon \}.
\end{eqnarray*}
Let
\[
\pi_{j}: M_{r_{j}} \to A''
\]
denote the canonical supporting $*$-homomorphism for $\tilde{\varphi}_{j}$ (cf.\ \ref{order-zero-structure}), 
so that we have 
\[
\tilde{\varphi}_{j}(x) = \tilde{\varphi}_{j}(\be_{M_{r_{j}}}) \pi_{j}(x)  \mbox{ for all } x \in M_{r_{j}}.
\]
We estimate that
\begin{eqnarray*}
\lefteqn{\| \tilde{\varphi} \tilde{\psi}(b) \tilde{\varphi} \tilde{\psi}(b') \| }\\
& \stackrel{\ref{product-domination}}{\ge} & \| \tilde{\varphi}_{j} \tilde{\psi}_{j}  (b) \tilde{\varphi}_{j} \tilde{\psi}_{j} (b') \|^{2} \\
& = & \| \tilde{\varphi}_{j}  (\be_{M_{r_{j}}} )^{2}  \pi_{j}  ( \tilde{\psi}_{j}(b)   \tilde{\psi}_{j}(b') )\|^{2} \\
&  \ge &  \|      \pi_{j}  ( \tilde{\psi}_{j}(b')   \tilde{\psi}_{j}(b) ) \\
&&  \tilde{\varphi}_{j}  (\be_{M_{r_{j}} })^{2}  \pi_{j} ( \tilde{\psi}_{j}(b)   \tilde{\psi}_{j}(b') )\|^{2} \\
& = &\| \tilde{\varphi}_{j}(\tilde{\psi}_{j} (b) \tilde{\psi}_{j} (b'))\|^{4} \\
& = & \|\tilde{\varphi}_{j}(\be_{M_{r_{j}} })\|^{4} \|\tilde{\psi}_{j} (b) \tilde{\psi}_{j}  (b')\|^{4}
\end{eqnarray*}
for all $j \in \{1, \ldots, s\}$ and normalized $b,b' \in A$. 

It follows that for each $j \in I$ there are $c,c' \in \mathcal{F}$ 
such that $\|cc'\| < \varepsilon$ and 
\begin{eqnarray*}
\lefteqn{\|\tilde{\varphi}\tilde{\psi}(c) \tilde{\varphi} \tilde{\psi} (c')\| }\\
& \ge & \|\tilde{\varphi}_{j}(\be_{M_{r_{j}}})\|^{4}  \varepsilon^{-\frac{1}{2}} \| \tilde{\varphi} \tilde{\psi} (c) \tilde{\varphi} \tilde{\psi} (c') \| ,
\end{eqnarray*}
whence
\[
\| \tilde{\varphi}_{j} (\be_{M_{r_{j}} })\| \le \varepsilon^{\frac{1}{8}}
\]
and 
\[
\textstyle
\left\| \sum_{j \in I} \tilde{\varphi}_{j} (\be_{M_{r_{j}} }) \right\| \le (n+1) \varepsilon^{\frac{1}{8}}.
\]
Set 
\[
F:= \bigoplus_{j \in \{1, \ldots,s\} \setminus I} M_{r_{j}}
\]
and denote the respective components of $\tilde{\varphi}$ and $\tilde{\psi}$ by $\varphi$ and $\psi$, respectively. Then, we have 
\begin{eqnarray*}
\| b - \varphi \psi (b)\| & \le & \| b - \tilde{\varphi} \tilde{\psi} (b)\| - \| \tilde{\varphi} \tilde{\psi} (b) - \varphi \psi (b)\| \\
& \le & \varepsilon + \textstyle \left\| \sum_{j \in I} \tilde{\varphi}_{j}(\be_{M_{r_{j}} }) \right\|  \\
& \le & \varepsilon + (n+1) \varepsilon^{\frac{1}{8}} \\
& < & \varepsilon^{\frac{1}{16}}
\end{eqnarray*} 
for $b \in \mathcal{F}$. 

Moreover, if $c,c' \in \mathcal{F}$ satisfy $\|cc'\| < \varepsilon$, then by the definition of $\psi$ and $I$, we have
\begin{eqnarray*}
\lefteqn{\| \psi (c) \psi (c') \|^{4}} \\
& = & \max_{j \notin I} \| \tilde{\psi}_{j} (c) \tilde{\psi}_{j} (c') \|^{4} \\
& < &  \varepsilon^{- \frac{1}{2}} \| \tilde{\varphi} \tilde{\psi}(c) \tilde{\varphi} \tilde{\psi} (c') \| \\
& < &  \varepsilon^{- \frac{1}{2}} ( \|cc'\| + 2 \varepsilon) \\
& \le & 3 \varepsilon^{\frac{1}{2}} \\
& < & \varepsilon^{\frac{1}{4}},
\end{eqnarray*}
so
\[
\|\psi(c) \psi(c')\| < \varepsilon^{\frac{1}{16}},
\]
as desired.

Now if $A$ is not necessarily separable, then the set 
\[
\Gamma:= \{B  \mid B \subset A \mbox{ is a separable $C^{*}$-subalgebra with } \dimnuc B \le \dimnuc A\}
\]
is directed with the order given by inclusion. Equip
\[
\Lambda:= \Gamma \times \mathbb{N}
\]
with the alphabetical order, then $\Lambda$ is directed as well. Use the first part of the proof to obtain an almost order zero, piecewise contractive, $n$-decomposable   system of c.p.\ approximations
\[
(F_{B,\nu},\psi_{B,\nu},\varphi_{B,\nu})_{\nu \in \mathbb{N}}
\]
for each $B \in \Gamma$. Using Proposition~\ref{separable-wdr}, it is straightforward to check that this yields  an almost order zero, piecewise contractive, $n$-decomposable   system of c.p.\ approximations
\[
(F_{B,\nu},\psi_{B,\nu},\varphi_{B,\nu})_{(B,\nu) \in \Lambda}
\]
for $A$ as desired. 
\end{proof}

\begin{nots}
We shall call $(F_{\lambda},\psi_{\lambda},\varphi_{\lambda})_{\lambda \in \Lambda}$ as in Proposition~\ref{almost-order-zero-approximation} a system of almost order zero, piecewise contractive,  $n$-decomposable c.p.\ approximations.
\end{nots}

The next result says that, if $A$ is sufficiently noncommutative, then so may be chosen the piecewise contractive, $n$-decomposable c.p.\ approximations. This will be particularly useful in Section~\ref{dichotomy}, where we derive a dichotomy result for $C^{*}$-algebras with finite nuclear dimension.

\begin{prop}
\label{nonelementary-approximation}
Let $A$ be a separable $C^{*}$-algebra with $\dimnuc A \le n < \infty$, and let $k \in \mathbb{N}$ be given. Suppose that $A$ has no irreducible representation of rank strictly less than $k$.

Then, there is a system $(E_{\nu},\varrho_{\nu},\sigma_{\nu})_{\nu \in \mathbb{N}}$ of almost order zero, piecewise contractive, $n$-decomposable c.p.\ approximations of $A$ such that the irreducible representations of each $E_{\nu}$ have rank at least $k$. 
\end{prop}

\begin{proof}
Choose a system
\[
(\bar{E}_{\nu}, \bar{\varrho}_{\nu}, \bar{\sigma}_{\nu})_{\nu \in \mathbb{N}}
\]
of almost order zero, piecewise contractive, $n$-decomposable c.p.\ approximations of $A$. For each $\nu$, write 
\[
\bar{E}_{\nu} = E_{\nu} \oplus \check{E}_{\nu},
\]
where $\check{E}_{\nu}$ consists precisely of those matrix blocks of $\bar{E}_{\nu}$ with rank at most $k-1$. Let $\varrho_{\nu}$, $\check{\varrho}_{\nu}$, $\sigma_{\nu}$ and $\check{\sigma}_{\nu}$ denote the respective components of $\bar{\varrho}_{\nu}$ and $\bar{\sigma}_{\nu}$. 

Let $h \in A$ be a normalized strictly positive element, and set 
\[
\mu:= \limsup_{\nu \in \mathbb{N}} \|\check{\varrho}_{\nu}(h) \| = \| \check{\varrho}(h)\| ,
\]
where 
\[
\check{\varrho}:A \to \prod \check{E}_{\nu} / \bigoplus \check{E}_{\nu}
\]
is the c.p.c.\ order zero map induced by the $\check{\varrho}_{\nu}$. Using a free ultrafilter on $\mathbb{N}$ and the fact that $\prod \check{E}_{\nu}$ is $(k-1)$-subhomogeneous, it is straightforward to construct an irreducible representation
\[
\pi: \prod \check{E}_{\nu}/ \bigoplus \check{E}_{\nu} \to M_{l}
\]
for some $l \le k-1$ such that 
\[
\|\pi \check{\varrho}(h) \|= \mu.
\]
Since $\pi$ is a $*$-homomorphism, $\pi \check{\varrho}$ again is a c.p.c.\ order zero map, so by Theorem~\ref{order-zero-structure} there are a $*$-homomorphism
\[
\sigma: A \to M_{l}
\]
and $0 \le d  \le \be_{l} \in M_{l}$ such that 
\[
d \sigma(a) = \sigma(a) d = \pi \check{\varrho} (a)
\]
for any $a \in A$. But by our assumption on $A$, $\sigma$ has to be zero, whence 
\[
\| \check{\varrho}(h)\| = \mu= \| \pi \check{\varrho}(h)\| =0.
\]
Using that $\check{\varrho}$ is a positive map and that $h$ is a strictly positive element, it is straightforward to conclude that $\check{\varrho}=0$. It follows that $(E_{\nu},\varrho_{\nu},\sigma_{\nu})_{\nu \in \mathbb{N}}$ is a system of c.p.\ approximations with the right properties.
\end{proof}

\section{Kirchberg's covering number}
\label{cov}

\noindent
In \cite[Definition~3.1]{Kir}, Kirchberg introduced a new integer valued invariant for a unital $C^{*}$-algebra. This covering number is closely related to both decomposition rank and nuclear dimension. It does not directly generalize topological covering dimension though, since it measures how many order zero maps one needs to \emph{cover} a noncommutative space, as opposed to \emph{approximating} it. In this section we recall the definition and some facts from \cite{Kir}, and then compare the covering number to nuclear dimension.  

\begin{defin}
Let $A$ be a unital $C^{*}$-algebra and $n \in \mathbb{N}$. $A$ has covering number at most $n$, $\mathrm{cov} A \le n$, if the following holds: 

For any $k \in \mathbb{N}$, there are a finite-dimensional $C^{*}$-algebra $F$, $d^{(1)}, \ldots,d^{(n)} \in A$  and a c.p.\ map $\varphi:F \to A$ such that
\begin{enumerate}
\item $F$ has no irreducible representation of rank less than $k$
\item $\varphi$ is $(n-1)$-decomposable with respect to $F = F^{(1)} \oplus \ldots \oplus F^{(n)}$
\item $\be_{A} = \sum_{j=1}^{n} (d^{(j)})^{*}\varphi^{(j)}(\be_{F^{(j)}})d^{(j)}$.
\end{enumerate} 
\end{defin}

We recall some more facts and notation from \cite[Section~1]{Kir}.

\begin{nots}
If $A$ is a $C^{*}$-algebra and $\omega \in \beta\mathbb{N} \setminus \mathbb{N}$ a free ultrafilter, we denote by $A_{\omega}$ the ultrapower $C^{*}$-algebra
\[
A_{\omega}:= \ell_{\infty}(A) / c_{\omega}(A);
\]
we will often consider $A$ as a subalgebra of $A_{\omega}$ via the canonical embedding as constant sequences. We denote the two-sided annihilator of $A$ in $A_{\omega} \cap A'$ by $\mathrm{Ann}(A)$, i.e., 
\[
\mathrm{Ann}(A):= \{b \in A_{\omega} \mid b A = A b = \{0\}\}.
\]
Then,  $\mathrm{Ann}(A)$ is a closed ideal in $A_{\omega}\cap A'$; if $A$ is $\sigma$-unital, then $A_{\omega}\cap A'/\mathrm{Ann}(A)$ is a unital $C^{*}$-algebra, cf.\ \cite[Proposition~1.9]{Kir}.
\end{nots}

We shall see below that $\mathrm{cov}\, A \le \dimnuc A + 1$ for any sufficiently noncommutative unital $C^{*}$-algebra. However,  the results of \cite{Kir} show that the covering number of the quotient $A_{\omega}\cap A'/\mathrm{Ann}(A)$ often is much more relevant than that of $A$. The next result relates the nuclear dimension of $A$ to the covering number of $A_{\omega}\cap A'/\mathrm{Ann}(A)$. This will be particularly useful in Section~\ref{dichotomy}. It will also play a key role in \cite{NgWin:cfp}, where Ng and the first named author will show that finite nuclear dimension implies the corona factorization property, at least for sufficiently noncommutative unital $C^{*}$-algebras.

\begin{prop}
\label{cov-nonelementary-hereditary}
Let $A$ be a separable $C^{*}$-algebra with $\dimnuc A \le n < \infty$, and suppose that no hereditary $C^{*}$-subalgebra of $A$ has a finite-dimensional irreducible representation. Then, 
\[
\mathrm{cov}(A_{\omega} \cap A'/\mathrm{Ann}(A)) \le (n+1)^{2}.
\]
\end{prop}

\begin{proof}
By \cite[Proposition~1.9]{Kir}, $A_{\omega} \cap A'/\mathrm{Ann}(A)$ is unital. Lift the unit $\be$ to a positive normalized element $e \in A_{\infty} \cap A'$; $e$ may be represented by an approximate unit $(e_{\lambda})_{\lambda \in \mathbb{N}}$ of $A$. 

By Proposition~\ref{almost-order-zero-approximation} there is a system 
\[
(F_{\lambda}= F_{\lambda}^{(0)}\oplus \ldots \oplus F_{\lambda}^{(n)}, \psi_{\lambda}, \varphi_{\lambda})_{\lambda \in \mathbb{N}}
\]
of almost order zero, piecewise contractive, $n$-decomposable c.p.\ approximations for $A$. By passing to a subsequence of the approximations, and by rescaling, if necessary, we may assume that 
\[
\|\varphi_{\lambda} \psi_{\lambda}(e_{\lambda})\| \le 1 \; \forall \, \lambda \in \mathbb{N},
\]
that
\begin{equation}
\label{ww3c}
\varphi_{\lambda} \psi_{\lambda}(e_{\lambda}^{\frac{1}{2}}ae_{\lambda}^{\frac{1}{2}}) \to a \; \forall \, a \in A
\end{equation}
and that
\begin{equation}
\label{ww3b}
\|\varphi_{\lambda} \psi_{\lambda}(e_{\lambda}) - e_{\lambda}\| \to 0.
\end{equation}
Define c.p.c.\ maps 
\[
\tilde{\psi}_{\lambda}: A^{+} \to F_{\lambda}
\]
by
\[
\tilde{\psi}_{\lambda}(\, . \,):= \psi_{\lambda}(e_{\lambda}^{\frac{1}{2}}\, . \,e_{\lambda}^{\frac{1}{2}}).
\]
For each $\lambda$, we define
\[
\hat{\psi}_{\lambda}( \, . \,) := \psi_{\lambda}(e_{\lambda})^{-\frac{1}{2}} \tilde{\psi}_{\lambda}(\, . \,) \psi_{\lambda}(e_{\lambda})^{-\frac{1}{2}},
\]   
where the inverses are taken in the hereditary subalgebras $\tilde{F}_{\lambda}$ generated by the $\psi_{\lambda}(e_{\lambda})$, and
\[
\hat{\varphi}_{\lambda}( \, . \,) := \varphi(\psi_{\lambda}(e_{\lambda})^{\frac{1}{2}} \, . \, \psi_{\lambda}(e_{\lambda})^{\frac{1}{2}}),
\]   
then
\[
\hat{\varphi}_{\lambda} \hat{\psi}_{\lambda} = \varphi_{\lambda} \tilde{\psi}_{\lambda};
\]
moreover, the 
\[
\hat{\psi}_{\lambda}: A^{+} \to \tilde{F}_{\lambda}
\] 
are unital c.p.\ and the 
\[
\hat{\varphi}_{\lambda}:\tilde{F}_{\lambda} \to A
\]
are c.p.c.\ maps.

From \cite[Lemma~3.5]{KW}, we see that for $i \in \{0, \ldots,n\}$ and $\lambda \in \mathbb{N}$ and any  projection $p_{\lambda}  \in \tilde{F}_{\lambda}$, 
\[
\|\hat{\varphi}_{\lambda}(p_{\lambda} \hat{\psi}_{\lambda}(a)) - \hat{\varphi}_{\lambda}(p_{\lambda}) \hat{\varphi}_{\lambda} \hat{\psi}_{\lambda} (a)\| \le 3 \cdot \max \{\| \hat{\varphi}_{\lambda} \hat{\psi}_{\lambda} (a) - a\|, \, \| \hat{\varphi}_{\lambda} \hat{\psi}_{\lambda} (a^{2}) - a^{2}\| \},
\]
from which follows that 
\begin{equation}
\label{ww1}
\| \varphi^{(i)}_{\lambda} \tilde{\psi}^{(i)}_{\lambda} (a) - \varphi^{(i)}_{\lambda} \tilde{\psi}^{(i)}_{\lambda}(\be_{A}) \varphi_{\lambda} \tilde{\psi}_{\lambda}(a) \| \stackrel{\lambda \to \infty}{\longrightarrow} 0
\end{equation}
for any $a \in A$. 

Let $\tilde{F}_{\lambda,l}^{(i)}, \,  l \in \{1, \ldots,r_{\lambda}^{(i)}\}$ denote the matrix blocks of $\tilde{F}^{(i)}_{\lambda}$, and denote the components of $\varphi_{\lambda}^{(i)}$ and $\tilde{\psi}_{\lambda}^{(i)}$ by $\varphi^{(i)}_{\lambda,l}$ and $\tilde{\psi}^{(i)}_{\lambda,l}$ accordingly. 

By Proposition~\ref{nonelementary-approximation} and our hypotheses on $A$, for each $\lambda \in \mathbb{N}$, $i \in \{0, \ldots,n\}$ and $l \in \{1, \ldots,r_{\lambda}^{(i)}\}$ there is 
\begin{equation}
\label{ww3a}
(E^{(i)}_{\lambda,l,\nu}= E^{(i,0)}_{\lambda,l,\nu} \oplus \ldots \oplus E^{(i,n)}_{\lambda,l,\nu}, \varrho^{(i)}_{\lambda,l,\nu}, \sigma^{(i)}_{\lambda,l,\nu})_{\nu \in \mathbb{N}},
\end{equation} 
an almost order zero, piecewise contractive, $n$-decomposable system of c.p.\ approximations of $\mathrm{her}(\varphi^{(i)}_{\lambda,l}(e_{11})) \subset A$ with the additional property that the matrix blocks of each $E^{(i)}_{\lambda,l,\nu}$ have rank at least $k^{2}$. A moment's thought shows that there is a unital $*$-homomorphism 
\[
\theta^{(i)}_{\lambda,l,\nu}: M_{k} \oplus M_{k+1} \to E^{(i)}_{\lambda,l,\nu}
\]
for any $i, \lambda, l,\nu$. 

Let 
\[
\bar{\sigma}^{(i)}_{\lambda,l,\nu}: E^{(i)}_{\lambda,l,\nu} \to   M_{r^{(i)}_{\lambda,l}} \otimes   \mathrm{her}(\varphi^{(i)}_{\lambda,l}(e_{11}))  \cong \mathrm{her}(\varphi^{(i)}_{\lambda,l}(\be_{M_{r^{(i)}_{\lambda,l}}})) \subset A
\]
be the amplification of $\sigma^{(i)}_{\lambda,l,\nu}$, using the canonical supporting $*$-homomorphism $\pi^{(i)}_{\lambda,l}$ of $\varphi^{(i)}_{\lambda,l}$, i.e., 
\[
\bar{\sigma}^{(i)}_{\lambda,l,\nu}(e) := \sum_{s=1}^{r^{(i)}_{\lambda,l}} \pi^{(i)}_{\lambda,l}(e_{s1}) \sigma^{(i)}_{\lambda,l,\nu}(e) \pi^{(i)}_{\lambda,l}(e_{1s}) \mbox{ for } e \in E^{(i)}_{\lambda,l,\nu}.
\]
Note that 
\begin{equation}
\label{ww2}
[\bar{\sigma}^{(i)}_{\lambda,l,\nu}(E_{\lambda,l,\nu}^{(i)}), \varphi^{(i)}_{\lambda,l}(M_{r^{(i)}_{\lambda,l}})] = 0
\end{equation}
and that $\bar{\sigma}^{(i)}_{\lambda,l,\nu}$ is decomposable into a sum of $n+1$ c.p.c.\ order zero maps $\bar{\sigma}^{(i,j)}_{\lambda,l,\nu}$ with respect to $E^{(i)}_{\lambda,l,\nu}= \bigoplus_{j=0}^{n} E^{(i,j)}_{\lambda,l,\nu}$. 

Let us fix a finite subset $\mathcal{F} \subset A$ of positive normalized elements and $\varepsilon>0$. By \eqref{ww1}, \eqref{ww3b} and \eqref{ww3c}, we can find $\lambda_{0} \in \mathbb{N}$ such that, for all $\bar{\lambda} \ge \lambda_{0}$,  
\[
\| \varphi_{\bar{\lambda}} \tilde{\psi}_{\bar{\lambda}}(a) \varphi^{(i)}_{\lambda} (\be_{F^{(i)}_{\bar{\lambda}}}) - \varphi^{(i)}_{\bar{\lambda}} \tilde{\psi}^{(i)}_{\bar{\lambda}}(a) \| < \frac{\varepsilon}{4},
\]
\[
\|\varphi_{\bar{\lambda}} \psi_{\bar{\lambda}} (e_{\bar{\lambda}}) - e_{\bar{\lambda}}\| < \frac{\varepsilon}{2}
\]
and such that
\[
\|\varphi_{\bar{\lambda}} \tilde{\psi}_{\bar{\lambda}}(a) - a\| < \frac{\varepsilon}{4}
\]
for $a \in \mathcal{F}$. Choose some
\[
0< \zeta < \frac{1}{8(n+1)} \varepsilon.
\]
Fix some $\bar{\lambda} \ge \lambda_{0}$. By the choice of the approximations in \eqref{ww3a}, there is $\bar{\nu} \in \mathbb{N}$ such that
\[
\| \sigma^{(i)}_{\bar{\lambda},l,\bar{\nu}} \varrho^{(i)}_{\bar{\lambda},l,\bar{\nu}}(g_{\zeta,2 \zeta}(\varphi^{(i)}_{\bar{\lambda},l}(e_{11}))) - g_{\zeta,2 \zeta}(\varphi^{(i)}_{\bar{\lambda},l}(e_{11})) \| < \zeta
\]
for each $i \in \{0, \ldots,n\}$ and $l \in \{1, \ldots, r^{(i)}_{\bar{\lambda}}\}$. Here, we define $g_{\zeta,2\zeta} \in C([0,1])$  by
\[
g _{\zeta,2\zeta}(t):= 
\left\{ \begin{array}{ll}
0& \mbox{for } 0\le t \le \zeta, \\
1 & \mbox{for } t \ge 2\zeta, \\
\mathrm{linear} & \mathrm{ else}.
\end{array} \right.
\]
We then have
\[
\sum_{l} \bar{\sigma}^{(i)}_{\bar{\lambda},l,\bar{\nu}}(\be_{E^{(i)}_{\bar{\lambda},l,\bar{\nu}}}) \ge g_{\zeta,2\zeta}(\varphi^{(i)}_{\bar{\lambda}}(\be_{F^{(i)}_{\bar{\lambda}}})) - \zeta.
\]
For $i,j \in \{0, \ldots,n\}$ define
\[
E^{(i,j)}:= \bigoplus_{l} E^{(i,j)}_{\bar{\lambda},l,\bar{\nu}}
\]
and 
\[
\sigma^{(i,j)}:= \bigoplus_{l} \bar{\sigma}^{(i,j)}_{\bar{\lambda},l,\bar{\nu}};
\]
note that 
\[
\sigma^{(i,j)}: E^{(i,j)} \to A
\]
is a c.p.c.\ order zero map. Let $\theta^{(i,j)}$ denote the respective component of $\bigoplus_{l} \theta^{(i)}_{\bar{\lambda},l,\bar{\nu}}$. Define
\[
\bar{\Phi}^{(i,j)}: E^{(i,j)} \to A
\]
by
\[
\bar{\Phi}^{(i,j)}(x):= \sigma^{(i,j)}(x) \varphi^{(i)}_{\bar{\lambda}}\tilde{\psi}^{(i)}_{\bar{\lambda}}(\be_{A^{+}}) \mbox{ for } x \in E^{(i,j)}.
\]
Note that by \eqref{ww2}, $\bar{\Phi}^{(i,j)}$ is a c.p.c.\ order zero map, whence
\[ 
\Phi^{(i,j)}:= \bar{\Phi}^{(i,j)} \circ \theta^{(i,j)}
\]
also is a c.p.c.\ order zero map. We have
\begin{eqnarray*}
\lefteqn{\sum_{i,j=0}^{n} \Phi^{(i,j)}(\be_{M_{k}\oplus M_{k+1}}) }\\
& = & \sum_{i,j=0}^{n} \bar{\Phi}^{(i,j)}(\be_{E^{(i,j)}}) \\
& = & \sum_{i,j} \sum_{l} \bar{\sigma}^{(i,j)}_{\bar{\lambda},l,\bar{\nu}} ( \be_{E^{(i,j)}_{\bar{\lambda},l,\bar{\nu}}}) \varphi^{(i)}_{\bar{\lambda}} \tilde{\psi}^{(i)}_{\bar{\lambda}}(\be_{A^{+}}) \\
& = & \sum_{i} \sum_{l} \bar{\sigma}^{(i)}_{\bar{\lambda},l,\bar{\nu}} ( \be_{E^{(i)}_{\bar{\lambda},l,\bar{\nu}}}) \varphi^{(i)}_{\bar{\lambda}} \tilde{\psi}^{(i)}_{\bar{\lambda}}(\be_{A^{+}}) \\
& \ge & \sum_{i} (g_{\zeta,2\zeta}(\varphi^{(i)}_{\bar{\lambda}}(\be_{F^{(i)}_{\bar{\lambda}}})) \varphi^{(i)}_{\bar{\lambda}} \tilde{\psi}^{(i)}_{\bar{\lambda}}(\be_{A^{+}})- \zeta) \\
& \ge & \sum_{i} ( \varphi^{(i)}_{\bar{\lambda}}\tilde{\psi}^{(i)}_{\bar{\lambda}}(\be_{A^{+}}) - 2 \zeta) \\
& \ge & \varphi_{\bar{\lambda}} \tilde{\psi}_{\bar{\lambda}}(\be_{A^{+}}) - (n+1) 2 \zeta \\
& = &  \varphi_{\bar{\lambda}} \psi_{\bar{\lambda}}(e_{\bar{\lambda}}) - (n+1) 2 \zeta \\
& \ge & e_{\bar{\lambda}} - \frac{\varepsilon}{2} - (n+1) 2 \zeta \\
& \ge & e_{\bar{\lambda}} - \varepsilon.
\end{eqnarray*}
Furthermore, we estimate for $a \in \mathcal{F}$ and $x \in E^{(i,j)}$ that
\begin{eqnarray*}
\lefteqn{\| [ \bar{\Phi}^{(i,j)}(x),a]\|} \\
& \le & \| [ \bar{\Phi}^{(i,j)}(x), \varphi_{\bar{\lambda}} \tilde{\psi}_{\bar{\lambda}}(a)]\| + 2 \frac{\varepsilon}{4} \|x\| \\
& = & \| \sigma^{(i,j)}(x) \varphi^{(i)}_{\bar{\lambda}}(\be_{F^{(i)}_{\bar{\lambda}}}) \varphi_{\bar{\lambda}} \tilde{\psi}_{\bar{\lambda}}(a) - \varphi_{\bar{\lambda}} \tilde{\psi}_{\bar{\lambda}}(a) \varphi^{(i)}_{\bar{\lambda}}(\be_{F^{(i)}_{\bar{\lambda}}})  \sigma^{(i,j)}(x)  \| + \frac{\varepsilon}{2} \|x\| \\
& \le & \| \sigma^{(i,j)}(x)  \varphi^{(i)}_{\bar{\lambda}} \tilde{\psi}^{(i)}_{\bar{\lambda}}(a) - \varphi^{(i)}_{\bar{\lambda}} \tilde{\psi}_{\bar{\lambda}}(a)  \sigma^{(i,j)}(x)  \| + 2 \frac{\varepsilon}{2} \|x\| \\
& \stackrel{\eqref{ww2}}{=} & \varepsilon \|x\|,
\end{eqnarray*}
from which follows that
\[
\|[\Phi^{(i,j)}(y),a]\| \le \varepsilon \|y\| \mbox{ for } y \in M_{k} \oplus M_{k+1}.
\]
Since $\mathcal{F}$ and $\varepsilon>0$ were arbitrary, and since the construction above works for any $\bar{\lambda} \ge \lambda_{0}$, it is now straightforward to construct c.p.c.\ order zero maps 
\[
\tilde{\Phi}^{(i,j)}: M_{k} \oplus M_{k+1} \to A_{\infty} \cap A'
\]
for $i,j= 0, \ldots,n$, satisfying
\[
\sum_{i,j} \tilde{\Phi}^{(i,j)}(\be_{M_{k}\otimes M_{k+1}}) \ge e.
\]
The $\tilde{\Phi}^{(i,j)}$ drop to c.p.c.\ order zero maps 
\[
\hat{\Phi}^{(i,j)}: M_{k} \oplus M_{k+1} \to A_{\omega} \cap A'/\mathrm{Ann}(A)
\]
satisfying
\[
\sum_{i,j} \hat{\Phi}^{(i,j)} (\be_{M_{k}\otimes M_{k+1}}) \ge \be.
\]
It follows that 
\[
\mathrm{cov}(A_{\omega} \cap A'/\mathrm{Ann}(A)) \le (n+1)^{2}.
\]
\end{proof}

Combining the idea of \cite[Proposition~3.5]{Kir} with the use of Proposition~\ref{nonelementary-approximation} as in the preceding proof, one can also show the following generalization of \cite[Proposition~3.5]{Kir}. 

\begin{prop}
\label{cov-nonelementary-hereditary-2}
Let $A$ be a separable unital $C^{*}$-algebra with $\dimnuc A \le n < \infty$, and suppose that  $A$ has no finite-dimensional irreducible representation. Then, 
\[
\mathrm{cov}(A) \le n+1.
\]
\end{prop}

\section{A dichotomy result}
\label{dichotomy}

\noindent
In this section we will combine Proposition~\ref{nonelementary-approximation} above with  \cite[Proposition~3.7]{Kir} to prove a dichotomy result for $C^{*}$-algebras with finite nuclear dimension: They either have nontrivial quasitraces, or they are weakly purely infinite. This statement becomes particularly satisfactory in the simple case. We first need some background results on lower semicontinuous (l.s.c.) traces.

\begin{prop}
\label{extending-traces-J}
Let $A$ be a separable $C^{*}$-algebra and $J \lhd A$ a closed ideal. Suppose $\tau$ is a densely defined l.s.c.\ trace on $J$. Then, $\tau$ extends to a (not necessarily densely defined) l.s.c.\ trace on $A$. 
\end{prop}

\begin{proof}
Choose an increasing approximate unit $(e_{\nu})_{\nu \in \mathbb{N}}$ for $J$. Using that $\tau$ is densely defined, a standard modification shows that we may even assume that $\tau (e_{\nu}^{\frac{1}{2}})< \infty$ for all $\nu \in \mathbb{N}$. 

Since $\tau$ is a trace and the $e_{\nu}$ are increasing,  for any $a \in A_{+}$ we obtain an increasing sequence of positive numbers 
\[
(\tau(e_{\nu}^{\frac{1}{2}}a e_{\nu}^{\frac{1}{2}}))_{\nu} 
\]
(these are all finite since they are dominated by the numbers $\tau(e_{\nu}^{\frac{1}{2}}) \|a\|$). We may thus define
\[
\bar{\tau}(a):= \lim_{\nu} \tau(e_{\nu}^{\frac{1}{2}}a e_{\nu}^{\frac{1}{2}}) \mbox{ for } a \in A_{+}.
\]  
Then, $\bar{\tau}$ is a well-defined map from $A_{+}$ to $[0,\infty]$. It is clear by lower semicontinuity that $\bar{\tau}$ extends $\tau$, that it is l.s.c.\ and that 
\[
\bar{\tau}(s \cdot a + t\cdot b) = s \cdot \bar{\tau}(a) + t \cdot \bar{\tau}(b) 
\]
if $\bar{\tau}(a), \bar{\tau}(b) <\infty$ and $s,t \in \mathbb{R}_{+}$. It remains to check that
\[
\bar{\tau}(x^{*}x) = \bar{\tau}(xx^{*})
\]
for all $x \in A$. To this end, note that for $x \in A$, $\mu \in \mathbb{N}$ and $\varepsilon >0$, we may choose $\nu_{0}$ so large that 
\[
\|e_{\mu}^{\frac{1}{4}} x^{*} (\be - e_{\nu}) x e_{\mu}^{\frac{1}{4}}\| < \frac{\varepsilon}{\tau(e_{\mu}^{\frac{1}{2}})}
\]
for any $\nu \ge \nu_{0}$ (this is where we use that $J$ is an ideal in $A$). We then estimate
\begin{eqnarray*}
\lefteqn{\tau(e_{\mu}^{\frac{1}{2}}x^{*}x e_{\mu}^{\frac{1}{2}})}\\
& = & \tau(e_{\mu}^{\frac{1}{2}}x^{*} e_{\nu} x e_{\mu}^{\frac{1}{2}}) +  \tau(e_{\mu}^{\frac{1}{4}} e_{\mu}^{\frac{1}{4}}x^{*} (\be - e_{\nu}) x e_{\mu}^{\frac{1}{4}} e_{\mu}^{\frac{1}{4}}) \\
& \le & \tau(e_{\mu}^{\frac{1}{2}}x^{*} e_{\nu} x e_{\mu}^{\frac{1}{2}}) + \frac{\varepsilon}{\tau(e_{\mu}^{\frac{1}{2}})} \cdot \tau(e_{\mu}^{\frac{1}{2}}) \\
& = & \tau(e_{\nu}^{\frac{1}{2}}x e_{\mu} x^{*} e_{\nu}^{\frac{1}{2}}) + \varepsilon \\
& \le & \tau(e_{\nu}^{\frac{1}{2}}x x^{*} e_{\nu}^{\frac{1}{2}}) + \varepsilon \\
& \le & \bar{\tau}(xx^{*}) + \varepsilon.
\end{eqnarray*}
Since $\mu$ and $\varepsilon$ were arbitrary, it follows that $\bar{\tau}(x^{*}x) \le \bar{\tau}(xx^{*})$; since the argument is symmetric in $x$ and $x^{*}$, we see that  $\bar{\tau}(x^{*}x) = \bar{\tau}(xx^{*})$, as desired.
\end{proof}

\begin{cor}
\label{extending-traces-her}
Let $A$ be a separable $C^{*}$-algebra and $B \subset A$ a hereditary $C^{*}$-subalgebra. If $\tau$ is a bounded  nontrivial  trace on $B$, then there is a (possibly unbounded) nontrivial l.s.c.\ trace $\tau'$ on $A$. 
\end{cor}

\begin{proof}
Let $J \lhd A$ be the (closed) ideal generated by $B$. By Brown's Theorem (\cite[Theorem~2.8]{Br}), $B \otimes \mathcal{K} \cong J \otimes \mathcal{K}$, since $B$ is full in $J$. 

Let $\Tr$ denote the standard l.s.c.\ trace on $\mathcal{K}$, then $\tau \otimes \Tr$ yields a densely defined nontrivial l.s.c.\ trace on $J \otimes \mathcal{K}$. Let $\bar{\tau}$ denote the restriction to $J$; it is straightforward to check that $\bar{\tau}$ again is densely defined, l.s.c.\ and nontrivial. By Proposition~\ref{extending-traces-J}, $\bar{\tau}$ extends to a l.s.c.\ trace $\tau'$ on $A$; since $\bar{\tau}$ is nontrivial, so is $\tau'$.
\end{proof}

\begin{prop}
\label{traceless-nonelementary-hereditary}
Let $A$ be a separable $C^{*}$-algebra and suppose $A$ has no nontrivial l.s.c.\ trace. 

Then, no hereditary $C^{*}$-subalgebra of $A$ has a finite-dimensional irreducible representation.
\end{prop}

\begin{proof}
If $B \subset A$ was a hereditary $C^{*}$-subalgebra with a finited-dimensional irreducible representation, then $B$ also had a (necessarily nontrivial) tracial state. By Corollary~\ref{extending-traces-her}, this would yield a nontrivial l.s.c.\ trace on $A$, a contradiction. 
\end{proof}

We are now prepared to prove the main result of this section.

\begin{thm}
\label{dichotomy-thm}
Let $A$ be a separable  $C^{*}$-algebra with $\dimnuc A \le n < \infty$. 

If $A$ has no nontrivial l.s.c.\ 2-quasitrace, then $A$ is weakly purely infinite.

In particular, if $A$ is simple, it is either strongly purely infinite, hence absorbs the Cuntz algebra $\mathcal{O}_{\infty}$, or it is stably finite with at least one densely defined trace. 
\end{thm}

\begin{proof}
Suppose $A$ has no nontrivial l.s.c.\ 2-quasitrace. By Proposition~\ref{traceless-nonelementary-hereditary}, no hereditary $C^{*}$-subalgebra of $A$ has a finite-dimensional irreducible representation. By Proposition~\ref{cov-nonelementary-hereditary} this yields
\[
\mathrm{cov}(A_{\omega} \cap A'/\mathrm{Ann}(A)) \le (n+1)^{2} < \infty.
\]
By \cite[Proposition~3.7]{Kir}, this implies that $A$ is weakly purely infinite.

For the second statement, suppose $A$ is simple but not purely infinite. Then, $A$ is not weakly purely infinite by \cite{KirRor:pi}, so $A$ admits a nontrivial l.s.c.\ 2-quasitrace. Therefore, $A$ contains a nonzero hereditary $C^{*}$-subalgebra $B$ with a bounded 2-quasitrace, which is a trace by \cite{Haa:quasitraces} or \cite{Kir:quasitraces} since  $B$ is nuclear. But then $B \otimes \mathcal{K}$ has a densely defined trace $\tau$. By Brown's Theorem, $A$ is a hereditary subalgebra of $B \otimes \mathcal{K}$, and it is straightforward to check that  $\tau$ restricts to a (nonzero) densely defined trace on $A$. 

For the statement that a simple purely infinite $C^{*}$-algebra absorbs $\mathcal{O}_{\infty}$ see \cite{Kir:ICM} or \cite[Theorem~7.2.6]{Ror:encyc}.
\end{proof}

\begin{rem}
\label{rem-dichotomy}
As of this moment, we do not know whether in the preceding result a traceless $C^{*}$-algebra $A$ will even be strongly purely infinite; this would imply that $A$ is $\mathcal{O}_{\infty}$-stable. The problem is closely related to the question whether finite nuclear dimension implies $\mathcal{Z}$-stability for sufficiently noncommutative $C^{*}$-algebras, cf.\ Conjecture~\ref{regularity-conjecture}. In fact, Theorem~\ref{dichotomy-thm} may be regarded as encouraging evidence to this effect. 
\end{rem}

\section{Examples}
\label{examples}

\noindent
In this section we list a number of examples for which we can compute or at least give bounds of their nuclear dimension. We will exhibit more examples in the subsequent sections.

\begin{examp}
We have already seen that decomposition rank dominates nuclear dimension, and that the two theories agree in the zero-dimensional and in the commutative case, and for continuous trace  $C^{*}$-algebras. This makes most examples of \cite[Section~4]{KW} accessible to nuclear dimension as well. In particular, for irrational rotation algebras $A_{\theta}$, we have 
\[
\dimnuc A_{\theta} = \left\{
\begin{array}{ll}
1 & \mbox{ if } \theta \mbox{ is irrational}\\
2 & \mbox{ if } \theta \mbox{ is rational.}
\end{array}
\right.
\]
\end{examp}

\begin{examp}
In \cite{TomsWin:minhom} it will be shown that, if $\alpha$ is a minimal homeomorphism of an infinite,  compact, finite-dimensional, metrizable  space $X$, then 
\[
\dimnuc (C(X) \rtimes_{\alpha} \mathbb{Z}) \le 2 \dim X +1.
\]
Examples suggest that this is not the best possible estimate in general (see above), and that the nuclear dimension of the crossed product should be bounded by $\max\{1,\dim X\}$, at least in the minimal case. However, for applications it often only matters whether or not the topological dimension is finite. 

For the decomposition rank, the latter estimate, i.e.,  
\[
\dr  (C(X) \rtimes_{\alpha} \mathbb{Z}) \le \max\{1, \dim X\},
\]
is known in special cases, e.g.\  when the action $\alpha$ is a minimal diffeomorphism on a compact smooth manifold $X$. The known proofs of such results, however, require the full strength of the classification theory for stably finite nuclear $C^{*}$-algebras.  The result of \cite{TomsWin:minhom} has the advantage that its proof is much simpler, and more conceptual. In particular, it does not factor through classification theorems of any kind. 
\end{examp}

\begin{examp}
\label{toeplitz-example}
Being an extension of $C(S^{1})$ by the compacts, the Toeplitz algebra $\mathcal{T}$ has nuclear dimension at most 2 by Proposition~\ref{quotext}. As of this moment, we do not know whether the precise value is 1 or 2 (it is not 0, since $\mathcal{T}$ is not AF). Since the Toeplitz algebra is infinite, hence not quasidiagonal, its decomposition rank is infinite. This in particular shows that decomposition rank and nuclear dimension do not agree. 
\end{examp}

\begin{examp}
In \cite{Ror:finite-infinite}, R{\o}rdam constructed a simple, separable, unital, and nuclear $C^{*}$-algebra containing a finite and an infinite projection. This example does not have a nontrivial (quasi-)trace, nor is it purely infinite, so by Theorem~\ref{dichotomy-thm} it has infinite nuclear dimension. 
\end{examp}

\section{Kirchberg algebras}  
\label{kirchberg-algebras}

\noindent
Below we will establish that classifiable simple purely infinite $C^*$-algebras  
have finite nuclear dimension.
It suffices to prove this for classical Cuntz algebras and then use inductive 
limit representations of classifiable algebras.
To this end let us collect some standard notation  and background results. 

Fix $n \in \N$, $n\ge2$ and recall that the Cuntz-Toeplitz algebra $\Tt_n$ is 
the universal $C^*$-algebra generated by $n$ isometries $T_1, \ld , T_n$ 
subject to the relations $T_i^*T_j = \delta_{ij} \be$, whereas the Cuntz 
algebra $\Oo_n$ is the universal $C^*$-algebra generated by $n$ isometries 
$S_1, \ld , S_n$ subject to the relations $S_i^*S_j = \delta_{ij} \be$ and 
$\sum_{i=1}^n S_iS_i^*=\be$.
Let $I=\{1,\ld ,n \}$ and $W_n=\bigcup_{k=0}^{\infty} I^{k}$ be the set of 
multi-indices or words in the alphabet $I$. For $\mu = i_1 \ld i_k  \in W_n$ 
let $|\mu|=k$ be the length of the word $\mu$ and define $S_{\mu}=S_{i_1} \ld S_{i_k}$, similarly $T_{\mu}=T_{i_1} \ld T_{i_k}$.
Every element $x$ in the $*$-algebra generated by the $S_i$ (respectively 
$T_i$) has a  representation as a finite linear combination of the form 
$x = \sum_{\mu,\nu} \alpha_{\mu,\nu} S_{\mu}S_{\nu}^*$ (respectively 
$x = \sum_{\mu,\nu} \alpha_{\mu,\nu} T_{\mu}T_{\nu}^*$). 

The full Fock space is defined by 
\[
\Gamma (n)= \bigoplus_{l =0}^{\infty} H^{\ot^l}, 
\]
where $H$ is an $n$-dimensional 
Hilbert space and $H^0:=\C \Omega$.
Fixing an orthonormal basis $e_1,\ld ,e_n$ of $H$ gives the orthonormal basis 
$e_{\mu}=e_{i_1} \ot e_{i_2} \ot \ld \ot e_{i_k}$ of $\Gamma(n)$, where $\mu=i_1 \ld i_k$ runs 
through $W_n$. In fact we may as well define $\Gan = \ell^2(W_n)$.
We denote by $M_{\infty}$ the $*$-algebra spanned by the matrix units $e_{\mu,\nu}$, where $\mu,\nu \in W_n$. 
Clearly, $M_{\infty} \subs \overline{M}_{\infty} = \Kk (\Gamma (n))\subs B(\Gamma (n))$.

As is well-known (\cite{Ev}), $\Tt_n$ acts faithfully on $\Gan$ with generators 
$T_i \xi = e_i \ot \xi$. It contains the matrix units 
$e_{\mu,\nu}=T_{\mu}(\be -\sum_{i=1}^n T_iT_i^*)T_{\nu}^*$ and hence the ideal 
of compact operators giving the exact sequence
$$
0 \to \Kk \to \Tt_n \to \Oo_n \to 0.
$$
As in \cite{SZ} we can write 
$$
T_{\mu}T_{\nu}^*=\sum_{i=0}^{\infty} e_{\mu,\nu} \ot \be_{H^{\ot^i}}= \sum_{i=0}^{\infty} e_{\mu,\nu} \ot \be_{i},
$$
where the sum is to be taken in the strong topology. 
The map $$\Lambda (x) = \sum_{i=0}^{\infty} x \ot \be_{H^{\ot^i}}$$ defined 
for matrix units $x$ may be regarded as an unbounded completely positive 
map $$\Lambda : M_{\infty} \to \Tt_n.$$
For a fixed integer $k>0$ define the cut-off Fock space 
$$\Gakn :=  \bigoplus_{l =0}^{k-1} H^{\ot^l}.$$  
It gives rise to the factorization
$$
\Gan \cong \Gakn \ot \Gank
$$ 
via $e_{\mu} \leftrightarrow e_{u} \ot e_{\bar{\mu}}$, where
$\mu = u \bar{\mu}$ and $|\bar{\mu}|$ is the largest multiple of $k$ below or equal to $|\mu|$. 
Similarly, if $k_1 | k_2$ then $\Gamma_{k_2}(n) \cong \Gamma_{k_1}(n) \ot \Gamma_{k_2/k_1}(n^{k_1})$.

Corresponding to the first factorization above we consider the $C^*$-algebra 
$$\Aa_k:=B(\Gakn) \ot \Tt_{n^k}.$$ Since $\dim \Gakn = 1+n+ \ld + n^{k-1} = \frac{n^k-1}{n-1} =:d_k$, 
this algebra is just $M_{d_k}(\Tt_{n^k})$. As shown in \cite{Kr} it is also generated by 
periodic weighted shifts but we don't need this description here. 

Important for us is that $\Aa_k$ contains $\Tt_n$. Indeed, denoting the generators 
of $\Tt_{n^k}$ by $\hat{T}_{v}$, where $v \in W_n$ with $|v|=k$, the generators $T_1, \ldots, T_n$ of $\Tt_n$ 
have the following matrix representation in $\Aa_k$:
\begin{eqnarray*}
T_i
&=&
\sum_{j=0}^{\infty} e_{i,0} \ot \be_{j}  \;\;\;\;\;\; (\textup{in } \Gamma (n))\\
&=&
\left(\sum_{j=0}^{k-2} e_{i,0} \ot \be_{j} \right) \ot \be_{\Gamma (n^k)} + \sum_{|w|=k-1} 
e_{0,w} \ot \hat{T}_{iw}. \\
\end{eqnarray*}
Similarly, if $k_1|k_2$ then the generators of $\be \ot \Tt_{n^{k_1}} \subs \Aa_{k_1}$ 
lie in $\Aa_{k_2}$ and $B(\Gamma_{k_1}(n)) \ot 1 \subset B(\Gamma_{k_1}(n)) \ot B(\Gamma_{k_2/k_1}(n^{k_1})) \cong B(\Gamma_{k_2}(n))$ so that $\Aa_{k_1} \subs \Aa_{k_2}$. 
If $k_1 < k_2 < \ld$ is a sequence of positive integers such that $k_i|k_{i+1}$ then 
$\Aa((k_i))=\overline{\bigcup_i \Aa_{k_i}}$ is a subalgebra of $B(\Gan)$. 

%In any case $\Aa_k/\Kk \sim M_{d_{k}}(\Oo_{n^{k}})$

Now let $Q(\Gan)$ be the Calkin algebra $B(\Gan)/\Kk$ with quotient homomorphism 
$q: B(\Gan) \to Q(\Gan)$ so that $q(\Tt_n) = \Oo_n$ and $q(\Aa_k)=M_{d_k}(\Oo_{n^k})$. 
Notice that the quotient $q(\Aa((k_i)))=\Aa((k_i))/\Kk$ is an inductive limit 
$\Bb((k_i))=\lim_i M_{d_{k_i}}(\Oo_{n^{k_i}})$, which is a simple nuclear purely 
infinite $C^{*}$-algebra. 

Moreover, there is a canonical unital inclusion 
$\Oo_{n^k} \hookrightarrow \Oo_{n}$ given on generators by $s_{v} \mapsto s_v$, 
where $v \in W_n$ with $|v|=k$. (We think of $\Oo_{n^k}$ as being generated by 
the isometries $s_{v}=s_{i_1} \ldots s_{i_k}$.) We obtain a unital embedding 
$M_{d_k}(\Oo_{n^k}) \hookrightarrow M_{d_k}(\Oo_{n})$.

It is known from classification theory that a matrix algebra of the form $M_r(\Oo_s)$ 
is isomorphic to $\Oo_s$ if $r$ and 
$s-1$ are relatively prime (\cite{Ror:encyc} 8.4.11(i)). Since 
$$
d_k=1+n+n^2+ \ld + n^{k-1} \equiv 1+1+\ld+1 = k \mod (n-1)
$$
there are certainly infinitely many $k$ satisfying $M_{d_k}(\Oo_{n}) \cong \Oo_{n}$.

We will need the following variant of the unbounded completely positive map $\Lambda$. 
Define $\Lambda_k : M_{\infty} \to B(\Gan)$ by 
$$
\Lambda_k (x) = \sum_{l=0}^{\infty} x \ot \be_{H^{\ot^{kl}}}=\sum_{l=0}^{\infty} x \ot \be_{kl}.
$$
Clearly, $\Lambda=\Lambda_1$. 

\begin{lem} \label{lambdak}
In the notation above we have:
\begin{enumerate}
\item $\Lambda_k (M_{\infty}) \subs \Aa_k \cong M_{d_k}(\Tt_{n^k})$.
\item For a non-negative integer $r$ let 
$$
\Gamma_{r,r+k}:=\bigoplus_{l=r}^{r+k-1} H^{\ot^k}=  
\Gakn \ot H^{\ot^r},
$$ 
so that $\Gamma_{0,k}=\Gamma_k$ and $B(\Gamma_{r,r+k}) \cong M_{n^rd_k}$. 
Then, $\Lambda_k|_{B(\Gamma_{r,r+k})}$ is a $*$-homomorphism. 
\end{enumerate}
\end{lem}

\begin{proof}
(i) Given $\mu, \nu \in W_n$ there are unique decompositions $\mu=u\bar{\mu}$ and 
$\nu=v\bar{\nu}$ such that $|u|, |v| <k$ and $|\bar{\mu}|,|\bar{\nu}|$ are multiples of $k$. 
Then 
\begin{eqnarray*}
\Lambda_k(e_{\mu , \nu})
&=& 
\sum_{l=0}^{\infty} e_{\mu , \nu} \ot \be_{lk} \\
&=& 
e_{u,v} \ot \sum_{l=0}^{\infty} e_{\bar{\mu} , \bar{\nu}} \ot \be_{lk} \\
&=& 
e_{u,v} \ot \hat{T}_{\bar{\mu}}\hat{T}_{\bar{\nu}}^*
\end{eqnarray*}
which proves the claim. \\

\medskip

(ii) This follows because for $x \in B(\Gamma_{r,r+k})$ the summands 
$$x, x \ot \be_k, x \ot \be_{2k}, \ld$$ of $ \Lambda_k (x)$ act $*$-homomorphically on the 
pairwise orthogonal subspaces $$\Gamma_{r,r+k}, \Gamma_{r+k,r+2k}, \Gamma_{r+2k,r+3k}, \ld $$ 
respectively. 
\end{proof}

We denote the projection onto $\Gamma_{r,r+k}$ by $P_{r,r+k}$. Define $P_k=P_{k,2k}$ and 
$Q_k=P_{\lceil k/2 \rceil +k , \lceil k/2 \rceil +2k}$,
where, as usual, 
$ \lceil k/2 \rceil =\inf \{ n \in \Z \mid n \geq k/2 \}$. 

\medskip

We now define the following positive $k \times  k$ matrices. 
For $k$ even let $l=k/2$ and define:
$$
\kappa_k=[\kappa_{i,j}]=\frac{1}{l+1}
\left[
\begin{array}{ccccccccc}
1      & 1      & \ld &     &     &  &  \ld & 1 & 1\\
1      & 2      & \ld &     &     &  &      & 2 & 1\\
1      & 2      &     &     &     &  &      &   & \vdots  \\
\vdots & \vdots &     &     &     &  &      &   &  \\
       &        &     &  l  &  l  &  &      &   &  \\
       &        &     &  l  &  l  &  &      &   &  \\ 
       &        &     &     &     &  &      &   &  \\
\vdots &        &     &     &     &  &      &   & \vdots \\
1      &        &     &     &     &  &      &   & 1 \\
1      &  1     & \ld &     &     &  & \ld & 1 & 1 \\
\end{array}
\right].
$$
For $k$ odd let $l= \lceil k/2 \rceil$ and define:
$$
\kappa_k=[\kappa_{i,j}]=\frac{1}{l+1}
\left[
\begin{array}{ccccccccc}
1      & 1      & \ld &     &     &     &  \ld & 1 & 1\\
1      & 2      & \ld &     &     &     &      & 2 & 1\\
1      & 2      &     &     &     &     &      &   & \vdots  \\
\vdots & \vdots &     &     &     &     &      &   &  \\
       &        &     & l-1 & l-1 & l-1 &      &   &  \\
       &        &     & l-1 &  l  & l-1 &      &   &  \\ 
       &        &     & l-1 & l-1 & l-1 &      &   &  \\
       &        &     &     &     &  &      &   &  \\
\vdots &        &     &     &     &  &      &   & \vdots \\
1      &   2    &     &     &     &     &      &   & 1 \\
1      &  1     & \ld &     &     &     & \ld  & 1 & 1 \\
\end{array}
\right].
$$
Since square matrices with all entries equal to 1 are positive, it is easy to see that
the above matrices are positive contractions.

Regarding $x \in B(\Gamma_{r,r+k})$ as a $k \times k$ operator matrix 
$x=[x_{i,j}]$, where $x_{i,j} \in B(H^{\ot^{r+j-1}}, H^{\ot^{r+i-1}})$
we infer that the Schur multiplication $\kappa_k * [x_{i,j}] =[\kappa_{i,j} x_{i,j}]$ 
defines a completely positive contraction.

With this at hand we define the following completely positive maps
$$
\psi_k : \Tt_n \to  B(\Gamma_{k,2k}) \oplus B(\Gamma_{\lceil k/2 \rceil+k,\lceil k/2 \rceil +2k}) 
$$
and 
$$
\vphi_k :  B(\Gamma_{k,2k}) \oplus B(\Gamma_{ \lceil k/2 \rceil+k, \lceil k/2 \rceil +2k}) \to \Aa_k = M_{d_k} (\Tt_{n^k}) \subs B(\Gan)
$$
by
$$
\psi_k (x) = \kappa_k(P_k x P_k) \oplus \kappa_k(Q_k x Q_k)
$$
and
$$
\vphi_k (x \oplus y) = \Lambda_k (x) + \Lambda_k (y)
$$
Clearly $\|\psi_k\|=1$ and $\|  \vphi_k \| =2$.
Finally we consider the composition $q \circ \vphi_k \circ \psi_k$.

\begin{prop} \label{factor}
For $\mu , \nu \in W_n$ fixed we have
$$
q \circ \vphi_k \circ \psi_k (T_{\mu} T_{\nu}^*) \to s_{\mu} s_{\nu}^*,
$$
as $k \to \infty$, 
where $s_{\mu}=q(T_{\mu})$ are the generators of $\Oo_n$ in the Calkin algebra 
$Q(\Gan)$. 
\end{prop}

\begin{proof}
Define the $\N_0 \times \N_0$ matrices 
$$
A_k=0_k \oplus \kappa_k \oplus \kappa_k \oplus \ldots 
$$
and 
$$
B_k=0_k \oplus 0_l \oplus \kappa_k \oplus \kappa_k \oplus \ldots ,
$$
where $l= \lceil k/2 \rceil$ and $0_k$ and $0_l$ denote the $k \times k$ resp.\ $l \times l$ zero matrices.
One checks that the entries $\sigma_{i,j}$ of the matrix $A_k+B_k$ verify
$\sigma_{i,i}=1$ and $|\sigma_{i,i+p} - 1| \leq \frac{2+p}{l+1} \leq \frac{2(2+p)}{k}$, provided $i>k$ and $0<p<l$.
Regard every operator on $\Gamma(n)$ as an operator matrix $[x_{i,j}]$, where
$x_{i,j} \in B(H^{{\ot}^j},H^{{\ot}^i})$. Then further inspection shows that 
$$
\vphi_k \circ \psi_k (T_{\mu} T_{\nu}^*)= (A_k + B_k) * (T_{\mu} T_{\nu}^*),
$$  
where $*$ denotes again Schur multiplication. Thus provided $k$ is large compared to 
$|\mu|$ and $|\nu|$ 
we have 
$$
\vphi_k \circ \psi_k (T_{\mu} T_{\nu}^*) = \sum_{r=0}^{\infty} \sigma_{|\mu|+r, |\nu|+r} e_{\mu,\nu} \ot I_r.
$$
By passing to the Calkin algebra (i.e.\;applying $q$) we obtain  
$$
\|s_{\mu}s_{\nu}^* - q \circ \vphi_k \circ \psi_k (T_{\mu} T_{\nu}^*)\| \leq 2 (2+\big| |\mu|-|\nu| \big|) k^{-1}.
$$
Letting $k$ tend to infinity concludes the proof.
\end{proof}

\begin{rem}
Alternatively, following along the lines of the proof of Theorem~\ref{wdrasdim}, we could define 
$\psi_k(x)= h_0 P_k x P_k h_0 \oplus h_1 Q_k x Q_k h_1$, where $h_0$ and $h_1$ are suitable 
positive diagonal matrices and then argue as in \ref{wdrasdim}.
\end{rem}

\begin{thm} \label{On} 
We have $\dimnuc \Oo_n = 1$ for $n=2, 3,\ldots$ and $\dimnuc \Oo_{\infty} \leq 2$.
\end{thm}

\begin{proof}
Let $\{x_1, x_2, \ld ,x_N\}$ be a finite subset of $\Oo_n$ and $\eps >0$. We need to 
find a finite dimensional $C^*$-algebra of the form $F=F^{(0)} \oplus F^{(1)}$, a c.p.c.\ map  
$\psi : \Oo_n \to F$ and $\vphi : F \to \Oo_n$ c.p.\ such that 
$\vphi |_{F^{(0)}} $ and  $\vphi |_{F^{(1)}} $ are both order zero contractions and such that 
$\| x_i -\vphi \circ \psi (x_i) \| < \eps$ for $i=1,\ld,N$.

To begin the construction fix $ \rho : \Oo_n \to \Tt_n$, a u.c.p.\ lift of the quotient 
map $\Tt_n \to \Oo_n$, which exists by nuclearity of $\Oo_n$.
For suitable $k$ (to be determined shortly) let $ F=  B(\Gamma_{k,2k}) \oplus B(\Gamma_{[k/2]+k,[k/2]+2k})$ 
and define  $\psi = \psi_k \circ \rho : \Oo_n \to F$. 

Next observe that $q \circ \psi_k : F \to M_{d_k}(\Oo_{n^k})\subs q(B(\Gamma(n)))$ 
which we compose with the inclusion $ M_{d_k}(\Oo_{n^k}) \hookrightarrow M_{d_k}(\Oo_{n}) \cong \Oo_n$
(the latter for suitable $k$). Moreover, $M_{d_k}(\Oo_{n^k})$ contains the copy of $\Oo_n$ 
from the inclusion $\Tt_n \subs B(\Gan)$; we think of $\{x_1, x_2, \ld ,x_N\}$ as a subset 
in that copy and then know from \ref{factor} that 
we may find $k$ such that $\| x_i  - q \circ \vphi_k \circ \psi_k \circ \rho (x_i) \| < \eps/2$ 
for $i=1,\ld,N$ and such that $d_k$ and $n-1$ are relatively prime. (Note that $\vphi_k \circ \psi_k (C) \to 0$ as $ k \to \infty$ for any compact $C \in \Kk (\Gamma (n))$ so that the choice of $\rho$ does not really matter.)

Further, we may regard the inclusion given by 
$$
\Oo_n \hookrightarrow M_{d_k}(\Oo_{n^k})  \hookrightarrow M_{d_k}(\Oo_{n}) \cong \Oo_n
$$
as a unital $*$-endomorphism $\sigma$ of $\Oo_n$. It follows from classification theory that 
any such endomorphism  is approximately unitarily equivalent to the identity map on $\Oo_n$. 
Indeed, $\sigma$ is homotopic to $\textup{id}$ since it is implemented by a unitary $v$ in $\Oo_n$ in the sense that $\sigma(s_i)=vs_i$ for all $i=1,\ldots ,n$ 
and the unitary group of $\Oo_n$ is connected. By Kirchberg's Classification Theorem 
(\cite{Ror:encyc} 8.3.3(iii)) $\sigma$ and  $\textup{id}$ are asymptotically hence approximately unitarily equivalent.

Thus there is a unitary $u \in \Oo_n$ 
such that 
$$
\| u x_i u^* - \sigma (x_i) \|  < \eps/2.
$$
Define $\varphi (x) = u^*(\beta \circ \varphi_k(x))u$, where $\beta$ denotes the map from $M_{d_k}(\Tt_{n^k})$ to $\Oo_n$ discussed above. Then  $(F, \psi,  \vphi )$ is as desired.

\bigskip

The estimate for $\mathcal{O}_{\infty}$ follows since there is an obvious inductive limit representation 
$\Oo_{\infty} = \lim_{n \to \infty} \Tt_n$, and we know that $\dimnuc \Tt_n \leq 2$ 
because of the exact sequence $$ 0 \to \Kk \to \Tt_n \to \Oo_n \to 0$$ and \ref{quotext}(ii).
\end{proof}

Using Kirchberg--Phillips classification it can be shown that every 
Kirchberg algebra satisfying the UCT is an inductive limit of $C^*$-algebras 
of the form 
$$
(M_{k_{1}}\otimes \Oo_{n_1} \oplus \ldots \oplus M_{k_{r}} \otimes \Oo_{n_r}) \ot C(\T),
$$ 
where $n_i \in \{2,3,\ldots \} \cup \{\infty \}$ and $k_i \in \N$ (cf.\ \cite{Ror:encyc}, 8.4.11). 
Since the nuclear dimension of any such algebra is at most 5 by Proposition~\ref{permanence}, we obtain the following.

\begin{thm} \label{pissnu}
A Kirchberg algebra (i.e., a purely infinite, simple, separable, nuclear $C^*$-algebra) 
satisfying the UCT has nuclear dimension at most 5.
\end{thm}

\section{Roe algebras} \label{roe-algebras}

\noindent
In this section we explore a connection between the asymptotic dimension of 
a coarse space and the nuclear dimension of its uniform Roe algebra. 
Although both concepts may be defined for arbitrary coarse spaces 
(c.f.\;\cite{Roe}) we restrict ourselves to discrete metric spaces of bounded geometry, 
mostly for simplicity. 

Recall that a discrete metric space $(X,d)$ is said to be of bounded geometry if 
every ball $B_r(x)=\{ y \in X \mid d(x,y) \leq r \}$ of finite radius $r$ has finitely 
many elements, and the number of elements in all balls of a given radius is uniformly 
bounded, that is, $b_r:= \sup \{ | B_r(x)| \mid x \in X\} < \infty$ for all $r$. 
This class of coarse spaces includes many interesting examples, e.g.\ finitely 
generated discrete groups with a word length metric.

In this setting the uniform Roe algebra $UC_r^*(X)$ associated to $(X,d)$ can be 
defined as follows:

Consider complex matrices $[\alpha_{x,y}]$ indexed by $x, y \in X$ such that
\begin{enumerate}
\item there is $M \geq 0$ with $|\alpha_{x,y}|\leq M$ for all $x ,y \in X$ 
(i.e.\;$[\alpha_{x,y}]$ is uniformly bounded);
\item there is $r>0$ such that $\alpha_{x,y}=0$ whenever $d(x,y)>r$
(i.e.\;$[\alpha_{x,y}]$ has bounded width).
\end{enumerate}
The smallest $r$ in condition (ii) is called the width of the matrix
$a=[\alpha_{x,y}]$, denoted by $w(a)$. Any matrix satisfying (i) and (ii) 
defines a bounded operator on $\ell^2(X)$, again denoted by $a$. We have in 
fact the following elementary estimate. 

\begin{lem}\label{bound}
Let $a= [\alpha_{x,y}]$ be a matrix
satisfying (i) and (ii) above and let $$b(a):=b_{w(a)}= \sup \{ |B_{w(a)}(x)| \mid x \in X\}.$$ 
Then, $\|a\| \leq b(a) M$.
\end{lem}

\begin{proof}
For $(\beta_x) \in \ell^2(X)$ the sum $\gamma_x =\sum_y \alpha_{x,y}\beta_y$ 
is well-defined containing at most $b(a)$ many terms for each $x \in X$. 
Thus 
$$
|\gamma_x|^2  \leq b(a) M^2 \sum_{y \in B_{w(a)}(x)}|\beta_y |^2.
$$ 
Since $\sum_{y \in B_{w(a)}(x)}|\beta_y |^2 \leq b(a) \| (\beta_x)\|^2$
we obtain $\| (\gamma_x)\| \leq M b(a) \| (\beta_x)\|$.
\end{proof}

Define the Roe algebra $UC_r^*(X)$ of $(X,d)$ as the concrete $C^*$-algebra 
generated  by matrices satisfying (i) and (ii) above, that is, the closure 
of the set of such matrices. Note that if $a \in UC_r^*(X)$ has finite width 
then the matrix entries are uniformly bounded (by $\|a\|$). 

We next recall the definition of the asymptotic dimension of $(X,d)$.
Note first that by a uniform cover $\Uu$ of $X$ we mean a family of subsets of 
$X$ such that $\bigcup_{U \in \Uu} U =X$ and such that the diameters $d (U)$ of all $U \in \Uu$ 
are uniformly bounded. A cover $\Uu$ has multiplicity or order $n$
if there are $n+1$ different $U_0, \ldots , U_n \in \Uu$ such that 
$ U_0 \cap \ldots \cap U_n \neq \emptyset$ but any $n+2$ different elements 
in $\Uu$ have empty intersection.

\begin{defin}\label{asdymdef}
Let $(X,d)$ be a metric space. The asymptotic dimension $\asdim X$ does not 
exceed $n$ if for every uniform cover $\Uu$ there is a uniform cover $\Vv$ of order $n$ such that $\Uu$ refines $\Vv$ (i.e.\;every $U \in \Uu$ is contained in a $V \in \Vv$).
\end{defin}

A family $\Uu$ of subsets of $X$ is said to be $r$-discrete if the distance 
$d(U,U')>r$ for any two different $U,U'\in\Uu$. We need the following 
characterization of the asymptotic dimension which is part of \cite[Theorem~19]{BD}.

\begin{thm}\label{rdisjoint}
For a metric space $(X,d)$ the following conditions are equivalent.
\begin{enumerate}
\item The asymptotic dimension $\asdim X$ does not exceed $n$.
\item For arbitrarily large $r>0$ there exist $r$-discrete families 
$$\Uu^{(0)}, \ldots ,\Uu^{(n)}$$ 
of subsets of $X$ such that  $\Uu^{(0)}\cup  \ldots \cup \Uu^{(n)}$ is a uniform cover of $X$.
\end{enumerate}
\end{thm}

It is known that $UC_r^*(X)$ is nuclear if $(X,d)$ is a discrete metric space of 
bounded geometry and finite asymptotic dimension. We will reprove this by 
constructing explicit approximating nets in the proof of the main result of 
this section below. Notice that for $X=\Gamma$ a discrete group, its uniform Roe algebra $UC_r^*(\Gamma)$ is nuclear iff $\Gamma$ is exact. Also other approximation properties of $\Gamma$ can be formulated in terms of the uniform Roe algebra (\cite{Zac:transapp}).

\begin{lem}\label{prod}
Let $K$ be any index set and $(n_k)_{k \in K}$ a bounded family of positive integers.
Then $\prod_{k \in K}M_{n_k}$ is an AF algebra.
\end{lem} 

\begin{proof}
Without loss assume $(n_k)$ to be constantly equal to $n$ ($\prod_{k \in K}M_{n_k}$ is a finite direct sum of such). Then any partition 
$\Pp=\{P_1,\ldots,P_l\}$ of $K$ defines an embedding of 
$$\underbrace{M_n \oplus \ld \oplus M_n}_{l} \to \prod_{k \in K}M_{n_k}$$  
sending $x_1 \oplus \ldots \oplus x_l$ to the family constantly equal to $x_i$ on $P_i$
for $i=1 ,\ld ,l$. The union of all the ranges of these embeddings for all possible finite partitions is dense in $\prod_{k \in K}M_{n_k}$.
\end{proof}

\begin{thm} \label{wdrasdim}
Let $(X,d)$ be a discrete metric space of bounded geometry. 
Then $\dimnuc (UC_r^*(X)) \leq \textup{asdim}(X)$. 
\end{thm}

\begin{proof}
Let $r \in \N$ and choose, according to Theorem~\ref{rdisjoint}, 
uniform $r$-disjoint families $\Uu^{(0)}, \ldots ,\Uu^{(n)}$ 
such that $\bigcup_{i=0}^n \Uu^{(i)}$ covers $X$. We will define 
a completely positive contraction $$
\Psi_r: UC_r^*(X) \to A^{(0)} \oplus \ld \oplus A^{(n)},
$$ 
where $$A^{(i)}= \prod_{U \in \Uu^{(i)}} M_{|B_{r-1}(U)|}.$$ 
By Lemma \ref{prod} every $A^{(i)}$ is AF and moreover naturally contained 
in $UC_r^*(X)$. Let $$\Phi_r : A^{(0)} \oplus \ld \oplus A^{(n)} \to UC_r^*(X)$$ 
be defined by $$\Phi_r (a_0 \oplus \ldots \oplus a_n)= a_0 + \ldots + a_n.$$
Then $\Phi_r$ is a completely positive map which is $*$-homomorphic on every 
$A^{(i)}$. If we can show that 
$\Phi_r \circ \Psi_r (a) \to a$ for all $a \in  UC_r^*(X)$ then we are done 
since we can combine the  $\Phi_r$ and $\Psi_r$
with a standard approximating net $(\psi_{\lambda} , \phi_{\lambda})$ of 
$ A^{(0)} \oplus \ld \oplus A^{(n)}$, where the $ \phi_{\lambda} $
are order 0, in fact $*$-homomorphic using \ref{prod}.

In order to define $\Psi_r$ let
$$
h^{(i)}=\frac{1}{r}\sum_{U \in \Uu^{(i)}} \sum_{l=1}^r \chi_{B(U,l-1)},
$$
where $\chi_S$ denotes the characteristic function of $S$ and $$B(U,s)=\{ x \in X \mid d(x,U) \leq s \}.$$
Then $h^{(0)} , \ldots , h^{(n)}$ are commuting positive contractions; moreover
$$
\be \leq h:=\sum_{i=0}^n h^{(i)} \leq (n+1) \be.
$$
If $a \in UC^*_r(X)$ is given by the matrix 
$[\alpha_{x,y}]$ then $[h^{(i)} ,a]$ is given by the matrix $[ (h^{(i)}(x)-h^{(i)}(y))\alpha_{x,y}]$ 
and if $a$ has finite width $w(a)<r$ then this commutator has still the same width and by Lemma \ref{bound} 
it follows that 
\begin{eqnarray*}
\| [h^{(i)},a]\|
& \leq &
b(a) \sup \{ |h^{(i)}(x) - h^{(i)}(y) | \mid d(x,y) < w(a) \} \| a \| \\
& \leq &
\frac{w(a)}{r} b(a) \|a\| 
\end{eqnarray*}
and thus 
$$
\|[h,a]\| \leq \frac{n+1}{r} w(a) b(a) \| a \|.
$$
Now define $$h_i=\left(h^{(i)}h^{-1}\right)^{1/2}.$$
Since $[h^{-1} ,a]= h^{-1} [a,h]h^{-1}$ we have 
$$\| [ h^{-1},a] \| \leq \|h^{-1}\|^2 \| [h,a]\| \leq  \| [h,a]\|,$$ so that 
$$\|[h^{(i)}h^{-1} , a]\| \to 0$$ as $r \to \infty$. 

Approximating the function $t \mapsto t^{1/2}$ by polynomials and using 
$$\|[a,x^n]\| \leq n \|[a,x]\| \|x\|^{n-1}$$ for any $x$ we find that 
also $$\left\| \left[ \left(h^{(i)}h^{-1}\right)^{1/2}, a \right] \right\| = \| [h_i,a]\| \to 0$$ 
as $r \to \infty$, whenever $a \in  UC_r^*(X)$ has finite width. But since $\|h_i\|\leq 1$ 
it follows that this is true for all $a \in UC_r^*(X)$.

Now define the completely positive contraction
$$
\Psi_r (a) = h_0 a h_0 \oplus h_1 ah_1 \oplus \ldots \oplus h_n a h_n.
$$
Then 
$$
\Phi_r \circ \Psi_r (a) = \sum_{i=0}^n h_i a h_i. 
$$
Note that $\Phi_r \circ \Psi_r (\be)= \sum_{i=0}^n h_i^2 = \be$ so that $\Phi_r \circ \Psi_r$ is u.c.p., in particular a contraction.
Since for $a \in UC_r^*(X)$ of finite width we have 
\begin{eqnarray*}
\left\| \Phi_r \circ \Psi_r (a) - a \right\|
&= &
\left\|\sum_{i=0}^n h_i ah_i - \sum_{i=0}^n h_i^2 a \right\| \\
&=& 
\left\|\sum_{i=0}^n h_i [a,h_i]  \right\| \leq \sum_{i=0}^n \|[a,h_i]\| \to 0, 
\end{eqnarray*} 
it follows again that $\|\Phi_r \circ \Psi_r (a) - a \| \to 0$ for all $a \in UC_r^*(X)$ because $\| \Phi_r \circ \Psi_r \| \leq 1$ for all
$r$.
\end{proof}

\section{Outlook. Open problems.}
\label{outlook}

\noindent
In this final section we list a number of open problems and possible future developments of the theory.\\

It follows trivially from the definitions that decomposition rank dominates nuclear dimension, and our (purely)  infinite examples show that the two theories do not agree in general. One might ask, however, whether infiniteness is the only obstruction.

\begin{question}
If $A$ is a  $C^{*}$-algebra with $\dimnuc A < \infty$, and if $A$ has a faithful trace, do we have $\dimnuc A = \dr A$? Do we at least have $\dr A < \infty$?
\end{question}

We have by now established upper and lower bounds for the nuclear dimension of a number of examples; while for  many applications it is enough to know whether the dimension is finite or infinite, it would nontheless be more satisfying to know the precise values, at least for the most important examples. The problem is that in general it is hard to find lower bounds for the nuclear dimension -- and in this respect our theory behaves just as many other (both commutative and noncommutative) notions of dimension.  

\begin{prob}
Determine the precise value of the nuclear dimension of the Toeplitz algebra, the Cuntz algebra $\mathcal{O}_{\infty}$, and, more generally, of Kirchberg algebras satisfying the UCT. Is the nuclear dimension of the latter determined by algebraic properties of their $K$-groups, such as torsion?
\end{prob}

A conjecture of Toms relates various regularity properties for separable, simple, finite, unital, and nuclear $C^{*}$-algebras. Our nuclear dimension enables us to put this conjecture into a broader context.  

\begin{con}
\label{regularity-conjecture}
For a separable, simple, unital, infinite dimensional and nuclear $C^{*}$-algebra $A$, the following are equivalent:
\begin{enumerate}
\item $A$ has finite nuclear dimension.
\item $A$ is $\mathcal{Z}$-stable.
\item $A$ has strict comparison of positive elements.
\item $A$ has almost unperforated Cuntz semigroup.
\end{enumerate}
\end{con}

As we have mentioned earlier, it will be shown in \cite{TomsWin:minhom} that crossed products of continuous functions on compact and finite dimensional spaces by the integers via minimal homeomorphisms have finite nuclear dimension. One might ask for similar results when the underlying $C^{*}$-algebra is noncommutative, or when the group is more complicated.

\begin{prob}
Find conditions on $A$, $G$ and $\alpha$, under which $\dimnuc (A \rtimes_{\alpha} G)$ is finite. 
\end{prob}

It is an open problem whether it is possible to recover a coarse metric space (up to coarse equivalence) from its uniform Roe algebra. To make at least some progress in this direction, one might ask for a converse to Theorem~\ref{wdrasdim}.

\begin{question}
Suppose $X$ is a discrete metric space of bounded geometry. Do we have $\dimnuc (UC^{*}_{r}(X)) = \asdim X$?
\end{question}

If the preceding question has a negative answer, we face another, perhaps even more interesting task:

\begin{prob}
Characterize all coarse metric spaces the Roe algebras of which have finite nuclear dimension. Describe regularity properties at the level of spaces (or groups) which are implied by finite nuclear dimension of the associated Roe algebras.   
\end{prob}

We have seen that quite different examples of noncommutative topological spaces are accessible to nuclear dimension; it is therefore natural to try to apply our theory to objects of a more geometric nature, such as Connes' spectral triples, cf.\ \cite{Connes:compact-metric} and \cite{Connes:NCG}. 

\begin{prob}
Find examples of spectral triples $(A,\pi,D)$ for which the nuclear dimension of $A$ can be related to summability properties of $D$. 
\end{prob}

%\begin{thebibliography}{Steg}
\bibliographystyle{amsplain}

\end{document}